\documentclass[a4paper, 12pt]{amsart}
\pdfoutput=1
\usepackage[utf8]{inputenc}
\usepackage[margin=1in]{geometry} 
\usepackage{amsmath,amsthm,amssymb,amsfonts,amscd,tikz,tikz-cd,tikz-3dplot,bbm,nicematrix}
\usepackage{thmtools}
\usepackage{thm-restate}
\usepackage{mathtools}
\usepackage{mathrsfs}
\usepackage{indentfirst}
\usepackage{enumitem}
\usepackage{graphicx}
\usepackage{subcaption}
\usepackage{float}
\usepackage{array}
\usepackage{hyperref}
\setlist[enumerate,1]{label={(\alph*)}}
\title{Reduced superschemes and the combinatorics of toric supervarieties}
\author{Eric Jankowski}
\date{}

\newcommand{\J}{\mathscr{J}}

\renewcommand{\O}{\mathscr{O}}
\renewcommand{\P}{\mathscr{P}}

\renewcommand{\AA}{\mathbb{A}}
\newcommand{\PP}{\mathbb{P}}

\newcommand{\ZZ}{\mathbb{Z}}

\newcommand{\RR}{\mathbb{R}}
\newcommand{\CC}{\mathbb{C}}
\newcommand{\kk}{\mathbb{C}}

\newcommand{\tab}{\hspace{0.6cm}}

\renewcommand{\epsilon}{\varepsilon}
\newcommand{\m}{\mathfrak{m}}
\newcommand{\p}{\mathfrak{p}}
\newcommand{\q}{\mathfrak{q}}
\newcommand{\uq}{\mathfrak{uq}}
\newcommand{\h}{\mathfrak{h}}
\newcommand{\g}{\mathfrak{g}}
\renewcommand{\t}{\mathfrak{t}}
\newcommand{\s}{\mathfrak{s}}
\renewcommand{\d}[1]{\frac{\partial}{\partial #1}}

\renewcommand{\tilde}{\widetilde}

\renewcommand{\r}{\text{right}}

\newcommand{\0}{{\ol{0}}}
\newcommand{\1}{{\ol{1}}}

\DeclareMathOperator{\Hom}{Hom}
\DeclareMathOperator{\Ext}{Ext}

\DeclareMathOperator{\Spec}{Spec}

\DeclareMathOperator{\Der}{Der}
\DeclareMathOperator{\minSpec}{minSpec}
\DeclareMathOperator{\zdiv}{zdiv}
\DeclareMathOperator{\nonzdiv}{nonzdiv}
\DeclareMathOperator{\Ann}{Ann}
\DeclareMathOperator{\Proj}{Proj}

\DeclareMathOperator{\codim}{codim}
\let\Im\undefined
\DeclareMathOperator{\Im}{Im}
\DeclareMathOperator{\tr}{tr}

\DeclareMathOperator{\Lie}{Lie}
\DeclareMathOperator{\Gr}{Gr}
\DeclareMathOperator{\QGr}{QGr}

\DeclareMathOperator{\Stab}{Stab}

\DeclareMathOperator{\Nil}{Nil}
\DeclareMathOperator{\cl}{cl}
\DeclareMathOperator{\Int}{int}
\DeclareMathOperator{\height}{ht}
\DeclareMathOperator{\depth}{depth}

\newtheorem{theorem}{Theorem}
\newtheorem{lemma}{Lemma}[section]
\newtheorem{corollary}[lemma]{Corollary}
\newtheorem{proposition}[lemma]{Proposition}

\theoremstyle{definition}
\newtheorem{example}[lemma]{Example}
\newtheorem{remark}[lemma]{Remark}
\newtheorem{definition}[lemma]{Definition}
\newtheorem{propdef}[lemma]{Proposition/Definition}

\newenvironment{usethmcounterof}[1]{%
  \renewcommand\thetheorem{\ref{#1}}\theorem}{\endtheorem\addtocounter{theorem}{-1}}

\newcommand{\ol}[1]{\overline{#1}}

\begin{document}

\begin{abstract}
We propose new definitions of integral, reduced, and normal superrings and superschemes to properly establish the notion of a supervariety. We generalize several results about classical reduced rings and varieties to the supergeometric setting, including an equivalence of categories between certain toric supervarieties and decorated polyhedral fans. These decorated fans are shown to encode important geometric information about the corresponding toric supervarieties. We then investigate some naturally-occurring toric supervarieties inside the isomeric supergrassmannian, which we show admits a nice description as a decorated polytope.
\end{abstract}


\maketitle

\section{Introduction}

Toric geometry provides a powerful framework for studying algebraic varieties with combinatorial structures. Classically, toric varieties serve as a rich family of examples whose algebraic, geometric, and topological properties are well understood, leveraging the deep correspondence between toric varieties and polyhedral fans. In the realm of supergeometry, where algebraic varieties are enhanced with additional fermionic coordinates, a natural question then arises: to what extent does this story generalize to the supergeometric setting?

Let $T$ be a \textit{supertorus} (or \textit{quasitoral supergroup}), a supergroup whose even part $T_\0$ is a central torus. Supertori arise as the Cartan subgroups of quasireductive supergroups \cite{SerganovaQred}, and their representations have been studied by Shibata \cite{Shibata} and Gorelik et al.\ \cite{GSS}.

If $T$ acts generically transitively on a normal supervariety $X$ (i.e.\ $X$ admits an open orbit), we say $X$ is a \textit{toric supervariety}. This paper aims to continue the investigation of toric supervarieties initiated by the author in \cite{Jankowski}, wherein only the case of one odd dimension is considered.

We first take care in sections 2 and 3 to develop a theory of commutative superalgebra and algebraic supervarieties which is consistent with the classical theory. Although terms like algebraic supervariety, integral superdomain, and reduced superring have been used before, we propose novel definitions and hope to persuade the reader that they are the best super-analogs of their classical counterparts. Our definition of a supervariety most closely aligns with Sherman’s in \cite{Sherman}, which inspired the theory of supervarieties developed in this paper.

The main content of this paper is section 4, in which we develop a systematic way of building supertori and toric supervarieties. It will be seen that classification of toric supervarieties in full generality is ultimately a wild problem, but we obtain several partial results in various level of generality. In particular, if we impose an enhancement of the usual $(R_1)$ condition of Serre (Proposition/Definition \ref{propdef:HR1}), then we can classify toric supervarieties by their decorated fans (Definition \ref{def:LETDecoratedFan}).

\begin{restatable}{theorem}{thm:HR1TSVandDecoratedFanBijection}\label{thm:HR1TSVandDecoratedFanBijection}
There is a bijective correspondence between $(HR_1)$ toric supervarieties and decorated fans, up to isomorphism.
\end{restatable}

We will see that morphisms of $(HR_1)$ toric supervarieties are difficult to describe, so a categorical version of this theorem remains elusive. However, we offer a partial rectification in the form of the following theorem about ``large-orbit" toric supervarieties, i.e.\ those whose orbits have odd codimension no greater than their even codimension. All large-orbit toric supervarieties admit the $(HR_1)$ property.

\begin{restatable}{theorem}{thm:largeOrbitEquivalenceOfCats}\label{thm:largeOrbitEquivalenceOfCats}
    The category of large-orbit toric supervarieties is equivalent to the category of large-orbit decorated fans.
\end{restatable}

We also discuss decorated polytopes in section 4.9, and we prove a smoothness criterion for $(HR_1)$ toric supervarieties in section 4.11.

Large-orbit toric supervarieties appear frequently in nature, and seem to be the most faithful super-analogs of toric varieties. Notably, their decorated fans (and decorated polytopes) are considerably easier to describe, and many statements, such as the smoothness criterion, simplify drastically.

Finally, in section 5, we investigate a family of large-orbit toric supervarieties in the isomeric supergrassmannian. We show that their decorated polytopes occur naturally as decorated polytopes of the image of a morphism of supermanifolds akin to a momentum map, generalizing results about torus strata in \cite{GGMS}.

\subsection{Acknowledgements}
The author is grateful to Vera Serganova and Alexander Sherman for many insightful comments and discussions about algebraic supergeometry. This material is based upon work supported by the National Science Foundation Graduate Research Fellowship Program under Grant No.\ 2146752.

\section{Commutative Superalgebra}
\subsection{Basics and conventions}
\begin{definition}
    A \textit{superring} is a $\ZZ/2\ZZ$-graded ring $A = A_\0 \oplus A_\1$.
\end{definition}
If $a \in A_{\ol{i}}$ is homogeneous, we write $|a| = \ol{i}$ for the \textit{parity} of $a$. For our purposes, all superrings will be assumed unital, associative, and \textit{commutative} in the sense that $ab = (-1)^{|a| \cdot |b|} ba$ for homogeneous $a,b \in A$.

All morphisms must preserve parity, and all ideals and modules are compatibly $\ZZ/2\ZZ$-graded. Maximal (respectively, prime) ideals in particular correspond bijectively to maximal (respectively, prime) ideals of $A_\0$, or equivalently $A/(A_\1)$, since prime ideals contain all nilpotents.

Localization behaves the same as in the ordinary situation, except that the multiplicatively closed subset $S$ that we invert must be contained in $A_\0$. That is, if $S$ is such a subset, $S^{-1} A = S^{-1} A_\0 \oplus S^{-1} A_\1$ where we treat $A_\1$ as an $A_\0$-module. To avoid degenerate cases, we will assume $0 \notin S$, so that localization preserves non-zero divisors.

We write $\zdiv(A)$ for the set of zero divisors of $A$, $\nonzdiv(A)$ for the non-zero divisors, $A^\times$ for the units, and $\Nil(A)$ for the nilpotents.

\begin{definition}
    The \textit{total superring of fractions} of $A$ is $K(A) := \nonzdiv(A)^{-1} A$. We denote by $i : A \to K(A)$ the natural inclusion, by $J_A$ the ideal $i^{-1}((K(A)_\1)) \subseteq A$, and by $\ol{A}$ the \textit{underlying ring} $A/J_A$. If $S$ is any subset of $A$, we write $\ol{S}$ for its image in $\ol{A}$.
\end{definition}

We record the following super version of the prime avoidance lemma, whose proof is completely identical to that of the standard version.
\begin{lemma}\label{PrimeAvoidance}
    Let $I$ be an ideal of a superring $A$ such that $I \subseteq \p_1 \cup ... \cup \p_n$ for some prime ideals $\p_i$. Then $I \subseteq \p_i$ for some $i$.
\end{lemma}

\subsection{Reduced superrings and integral superdomains}
We introduce some terminology which has seen conflicting usage in the literature. We believe the following definitions are most in line with the classical situation. The common theme here is that badly-behaved elements should all arise from the odd part of the total superring of fractions.

\begin{definition}
    Let $A$ be a superring.
    \begin{enumerate}
        \item $A$ is a \textit{superfield} if $J_A = A \backslash A^\times$
        \item $A$ is an \textit{integral superdomain} if $J_A = \zdiv(A)$
        \item $A$ is \textit{reduced} if $J_A = \Nil(A)$ and $\ol{\zdiv(A)} = \zdiv(\ol{A})$
    \end{enumerate}
\end{definition}

In Lemma \ref{SuperfieldIntegralReducedLemma} we will provide several equivalent characterizations of these types of superrings. We first prove a lemma that will allow us to conflate $\ol{A_\p}$ and $\ol{A}_{\ol{\p}}$ for prime ideals $\p$ of a reduced superring $A$.

\begin{lemma}\label{LocalizationPreservesReducedJALemma}
    Let $A$ be a reduced superring and $S \subseteq A_\0$ a multiplicatively closed subset.
    \begin{enumerate}
        \item $J_{S^{-1} A} = S^{-1} J_A$
        \item $\ol{S^{-1} A} = \ol{S}^{-1} \ol{A}$
    \end{enumerate}
\end{lemma}

\begin{proof}
    \begin{enumerate}
        \item Since $\Nil(S^{-1} A) = S^{-1}\Nil(A) = S^{-1} J_A$, it suffices to show that $\Nil(S^{-1} A) = J_{S^{-1} A}$. Let $\frac{a}{s} \in S^{-1} A$ be nilpotent, so there is $t \in S$ such that $ta^n=0$. Then $ta \in \Nil(A) = J_A$, so $ta = \frac{\eta}{b}$ for $\eta \in (A_\1)$ and $b \in \nonzdiv(A)$. Then $b \cdot \frac{a}{s} = \frac{\eta}{st}$ in $S^{-1} A$ for $\frac{\eta}{st} \in (S^{-1} A_\1)$. Hence $\frac{a}{s} \in J_{S^{-1}A}$, since localization preserves non-zero divisors. Therefore $\Nil(S^{-1}A) = J_{S^{-1} A}$ and we are finished.

        \item Follows immediately from part (a). \qedhere
    \end{enumerate}
\end{proof}

\begin{lemma}\label{SuperfieldIntegralReducedLemma}
    Let $A$ be a superring.
    \begin{enumerate}
        \item The following are equivalent: \begin{enumerate}[label=(\roman*)]
            \item $A$ is a superfield
            \item $\ol{A}$ is a field
            \item $J_A = (A_\1)$ is the set of non-units
            \item $(A_\1)$ is a maximal ideal
        \end{enumerate}
        
        \item The following are equivalent: \begin{enumerate}[label=(\roman*)]
            \item $A$ is an integral superdomain
            \item $\ol{A}$ is an integral domain and $\ol{\zdiv(A)} = \zdiv(\ol{A})$
            \item $K(A)$ is a superfield
            \item $A$ is reduced and $\Spec A$ is irreducible
        \end{enumerate}

        \item The following are equivalent: \begin{enumerate}[label=(\roman*)]
            \item $A$ is a reduced superring
            \item $\ol{A}$ is a reduced ring and $\ol{\zdiv(A)} = \zdiv(\ol{A})$
            \item $J_A = \Nil(A)$ and $\zdiv(A) = \bigcup_{\p \in \minSpec A} \p$
            \item $\ol{A}$ is a reduced ring and $\zdiv(A) = \bigcup_{\p \in \minSpec A} \p$
            \item $S^{-1} A$ is a reduced superring for any multiplicatively closed subset $S \subseteq A_\0$
            \item $A_\m$ is a reduced superring for any maximal ideal $\m$
        \end{enumerate}
    \end{enumerate}
\end{lemma}

\begin{proof}
    \begin{enumerate}
        \item The equivalence of (iii) and (iv) is immediate since every maximal ideal contains $(A_\1)$. These equivalent conditions imply $(A_\1) = J_A$, which then implies (i) and (ii), which are equivalent because every maximal ideal contains $J_A$.

        \tab Let us now prove (i) $\implies$ (iv). It suffices to show that $J_A = (A_\1)$. Let $\eta \in J_A$, so we may write $\eta = \frac{\eta'}{a}$ for $\eta' \in (A_\1)$ and $a \in \nonzdiv(A)$. But $a \notin J_A$ implies $A$ is a unit, so indeed $\eta \in (A_1)$ and we are finished.

        \item First assume (i), so that $K(A) = A_{J_A}$ has maximal ideal $K(A) J_A = (K(A)_\1)$ and (iii) follows. For the converse, note that a zero divisor $a \in A$ does not receive an inverse in $K(A)$, so $a \in J_A$. Thus, we have established (i) $\iff$ (iii).

        \tab Now observe that (ii) just says $\ol{\zdiv(A)} = 0$, so indeed $J_A = \zdiv(A)$, which is (i). Conversely, since nonzerodivisors form a multiplicatively closed subset, $A$ integral implies $\ol{A}$ integral and so we obtain (i) $\iff$ (ii). Finally, (ii) $\iff$ (iv) follows from part (c).

        \item That (i) $\iff$ (ii) and (iii) $\iff$ (iv) is immediate from the definitions.  To see the equivalence of (i) and (ii) with (iii) and (iv), it suffices to show that if $\ol{A}$ is reduced, then $\ol{\zdiv(A)} = \zdiv(\ol{A})$ if and only if $\zdiv(A) = \bigcup_{\p \in \minSpec A} \p$.

        For the forward implication, let $a \in \zdiv(A)$, so there is $b \in A$ such that $b \notin J_A$ and $ab \in J_A$. Then $ab \in \p$ for all primes $\p$, but $b \notin \p'$ for some minimal prime $\p'$ since $b \notin \Nil(A)$. Hence $a \in \p'$, so $\zdiv(A) \subseteq \bigcup_{\p \in \minSpec A} \p$. That minimal primes consist of zero divisors is standard. Conversely, if $\zdiv(A) = \bigcup_{\p \in \minSpec A} \p$, then $\ol{\zdiv(A)} = \ol{\bigcup_{\p \in \minSpec A} \p} = \bigcup_{\p \in \minSpec A} \ol{\p} = \zdiv(\ol{A})$ since $\ol{A}$ is reduced.

        We now show (iii) $\implies$ (v). That $\Nil(S^{-1} A) = J_{S^{-1}A}$ is proven in Lemma \ref{LocalizationPreservesReducedJALemma}. For $\zdiv(S^{-1}A) = \bigcup_{\p \in \minSpec S^{-1}A} \p$, first take $\frac{a}{s} \in \zdiv(S^{-1}A)$, so $a \in \zdiv(A)$. Then $a$ belongs to a minimal prime $\p$ of $A$. If $\p \cap S = \emptyset$, then $S^{-1} \p$ is a minimal prime of $S^{-1} A$ and we are finished. Otherwise we may assume that $a$ belongs to no such minimal prime. Let $\q$ be a prime that does contain $a$, so $\q$ contains a minimal prime that meets $S$. Then $S^{-1}\q = S^{-1} A$, so $s^{-1} a$ is a unit in $S^{-1}A$, a contradiction. Hence $\zdiv(S^{-1}A) = \bigcup_{\p \in \minSpec S^{-1}A} \p$, so indeed (iii) $\implies$ (v). 
        
        Since (v) $\implies$ (vi) is immediate, it remains to show (vi) $\implies$ (i). First let $a \in \Nil(A)$, define $I = \{d \in A \mid ad \in J_A\}$, and pick a maximal ideal $\m$. Then $\frac{a}{1} \in \Nil(A_\m) = J_{A_\m}$, so $\frac{b}{t} \cdot \frac{a}{1} = \frac{\alpha}{s}$ for $\frac{\alpha}{s} \in ((A_\m)_\1)$ and $\frac{b}{t} \in \nonzdiv(A_\m)$. Hence there is $u \notin \m$ such that $u(abs-t\alpha) = 0$ in $A$. Since $u,b,s \notin \m$ and $ut\alpha \in J_A$, it follows that $ubs \in I \backslash \m$, so $I$ is not contained in $\m$. Therefore $I$ is the unit ideal, so indeed $a \in J_A$.
        
        Now if $A$ is not reduced, there is a nonzero $\xi \in J_A$ such that $\ol{\Ann_A(\xi)}$ is not contained in $\zdiv(\ol{A})$. Let $a \in \Ann_A(\xi)$ be such a ``bad" element, so that $\ol{a} \notin \zdiv(\ol{A})$. Moreover, let $\m$ be a maximal ideal containing $\Ann_A(\xi)$. Then $\Ann_A(a) \subseteq J_A \subseteq \m$, so the relationship $a \xi =0$ remains nontrivial in $A_\m$. This contradicts reducedness of $A_\m$, so indeed $A$ must be reduced as well. \qedhere
    \end{enumerate}
\end{proof}

We now collect some results about reduced superrings in the general and Noetherian cases.

\begin{lemma}\label{ReducedLocalizationMinimalPrimeSuperfieldLemma}
    Let $A$ be a reduced superring. Then $A_\p$ is a superfield for all $\p \in \minSpec A$.
\end{lemma}

\begin{proof}
    If $\p$ is a minimal prime, then $A_\p$ is a reduced local superring with a unique prime ideal $\p A_\p = \Nil(A_\p) = J_{A_\p}$ and hence a superfield.
\end{proof}

\begin{lemma}\label{ReducedNoetherianLocalizationsLemma}
    Let $A$ be a Noetherian reduced superring and $S \subseteq A_\0$ a multiplicative subset.
    \begin{enumerate}
        \item $K(A) \cong \prod_{\p \in \minSpec A} A_\p$
        \item $S^{-1} K(A) = K(S^{-1}A)$
    \end{enumerate}
\end{lemma}

\begin{proof}
    \begin{enumerate}
        \item By the universal property of localization, there are natural maps $K(A) \to A_\p$, since minimal primes consist only of zero divisors. Hence we get a map $K(A) \to \prod_{\p \in \minSpec A} A_\p$. By the prime avoidance lemma \ref{PrimeAvoidance}, any nonminimal prime contains a non-zero divisor, so $\Spec K(A) = \minSpec A$, a finite discrete set, as a subspace of $\Spec A$. Passing to spectra, it follows that $K(A) \to \prod_{\p \in \minSpec A} A_\p$ is an isomorphism.

        \item We know that $S^{-1}A$ is reduced by Lemma \ref{SuperfieldIntegralReducedLemma}(c). Hence $K(S^{-1}A)$ is the product of $(S^{-1}A)_{S^{-1}\p}$ over minimal primes $\p$ of $A$ such that $\p \cap S = \varnothing$. If $S$ meets $\p$, then $S^{-1} (A_\p) = 0$ since $\p A_\p = ((A_\p)_\1)$ and hence consists of only nilpotents. Therefore $S^{-1} K(A) \cong \prod_{\p \in \minSpec S^{-1} A} A_\p \cong K(S^{-1}A)$. \qedhere
    \end{enumerate}
\end{proof}

\begin{lemma}\label{EquivalentReducedNoetherianLemma}
    Let $A$ be a Noetherian superring. The following are equivalent:
    \begin{enumerate}[label=(\roman*)]
        \item $A$ is reduced
        \item $A_\p$ is a superfield for all $\p \in \minSpec A$, and $A \to \prod_{\p \in \minSpec A} A_\p$ is injective
    \end{enumerate}
\end{lemma}
\begin{proof}
    We have (i) $\implies$ (ii) from Lemmas \ref{ReducedLocalizationMinimalPrimeSuperfieldLemma} and \ref{ReducedNoetherianLocalizationsLemma}.
    
    Now assume (ii), and let $s \in A \backslash \bigcup_{\p \in \minSpec A} \p$. Then the image of $s$ in $\prod_{\p \in \minSpec A} A_\p$ is invertible and hence a non-zero divisor in $A$, so indeed $\zdiv(A) = \bigcup_{\p \in \minSpec A} \p$. We now show $\Nil(A) = J_A$. Let $a \in \Nil(A)$. Then $\frac{a}{1} \in J_{A_\p} \cong (J_A)_\p$, i.e.\ there is $s \notin \p$ such that $as \in J_A$. Let $I = \{s \in A \mid as \in J_A\}$, so $I \not\subseteq \p$ for any minimal prime $\p$. Using the prime avoidance lemma \ref{PrimeAvoidance}, let $s \in I \backslash \bigcup_{\p \in \minSpec A} \p$. As before, $s \in \nonzdiv(A)$, so indeed $a \in J_A$.
\end{proof}

When we generalize the commutative superalgebra of this section to algebraic supergeometry in the next, we will often want our spaces to contain smooth points. Part (c) of the following definition provides a suitably general algebraic model for such spaces.

\begin{definition}\label{def:FRsuperring}
    Let $A$ be a superring.
    \begin{enumerate}
        \item $A$ is \textit{fermionically regular} (or \textit{FR}) if it is of the form $\bigwedge M$ for some commutative ring $R$ and $R$-module $M$. In this situation, we have $J_A = (A_\1)$, and we may assume $R \cong \ol{A}$ and $M \cong J_A/J_A^2$.
        \item $A$ is \textit{generically fermionically regular} (or \textit{GFR}) if $A_\p$ is a FR superring for all minimal prime ideals $\p$.
        \item $A$ is \textit{GFRR} if it is GFR and reduced (so that each $A_\p$ is a FR superfield for $\p$ minimal).
    \end{enumerate}
\end{definition}

\subsection{Serre's criteria}
Let $A$ be a superring. Following \cite{Schmitt}, we say an element $t \in A_\0$ is $A$-regular if the multiplication map by $t$ is injective but not bijective. Likewise, we say $\xi \in A_\1$ is $A$-regular if the cohomology of the multiplication map by $\xi$ is trivial. An $A$-regular sequence is a sequence of homogeneous elements $r_1, ..., r_k$ such that each $r_i$ is regular in $A/(r_1, ..., r_{i-1})$.

\begin{definition}
    Let $A$ be a Noetherian local superring with maximal ideal $\m$.
    \begin{enumerate}
        \item $A$ is \textit{regular} if $\m$ is generated by an $A$-regular sequence
        \item The \textit{depth} of $A$ is the minimum $i\geq0$ such that $\Ext^i(A/\m, A) \neq 0$.
    \end{enumerate}
\end{definition}

The usual theorem of Rees still holds in this setting by the same argument:
\begin{lemma}
    The depth of a Noetherian local superring $(A,\m)$ is the length of any maximal $A$-regular sequence of even elements in $\m$.
\end{lemma}

\begin{definition}\label{def:SerresCriteria}
    Let $A$ be a Noetherian superring. We say $A$ satisfies
    \begin{enumerate}
        \item $(R_k)$ if $A_\p$ is regular for all prime ideals $\p$ of height $\leq k$
        \item $(S_k)$ if $\depth(A_\p) \geq \min(k, \height(\p))$ for all prime ideals $\p$
    \end{enumerate}
\end{definition}

\begin{lemma}
    A superfield is regular if and only if it is FR.
\end{lemma}

The following proposition justifies our usage of GFRR superrings as super-versions of reduced rings, as well as our definitions of the conditions $(R_k)$ and $(S_k)$.

\begin{proposition}
    A Noetherian superring $A$ is GFRR if and only if it satisfies $(R_0)$ and $(S_1)$.
\end{proposition}

\begin{proof}
    If $A$ is GFRR, let $\p$ be a minimal prime. Then $A_\p$ is a FR superfield, so $(R_0)$ holds. Now let $\p$ be a prime ideal with $\height(\p) \geq 1$. For $(S_1)$, it suffices to show that there is a nonzerodivisor in $\p A_\p$. We know by Lemma \ref{SuperfieldIntegralReducedLemma} that $A_\p$ is reduced, so we have $\zdiv(A_\p) = \bigcup_{\q \in \minSpec(A_\p)} \q = \bigcup_{q \in \minSpec(A), \q \subseteq \p} \q A_\p$. If $\p A_\p$ is contained in this set, then $\p$ is one of these minimal primes by the prime avoidance lemma \ref{PrimeAvoidance}. But $\height(\p) \geq 1$, so indeed $\p A_\p$ contains a nonzerodivisor.

    Conversely, suppose $A$ satisfies $(R_0)$ and $(S_1)$. We will show that $A_\p$ is reduced for all primes $\p$. If $\p$ is a minimal prime, then by $(R_0)$ we know $A_\p$ is regular. But $\p A_\p$ is the unique prime ideal of $A_\p$, so $\p A_\p = \Nil(A_\p)$. Regularity implies that $\Nil(A_\p) = ((A_\p)_\1)$, so $A_\p$ is a superfield and therefore is reduced.
    
    Now if $\height(\p) \geq 1$, we know by $(S_1)$ that $\p A_\p$ contains a nonzerodivisor, i.e.\ an element $a$ for which the localization map $A_\p \to A_\p[a^{-1}]$ is injective. Since any superring, particularly $A_\p[a^{-1}]$, is contained in the product of its localizations at prime ideals, we know that $A_\p$ is a subring of a product of localizations of $A$ at primes $\q \subsetneq \p$ with $a \notin \q$. Since $\height(\q) < \height(\p)$, we may induct on the height to find that $A_\p$ is a subring of a product of reduced superrings which are localizations of $A_\p$ at minimal primes. It follows by Lemma \ref{EquivalentReducedNoetherianLemma} that $A_\p$ (and hence also $A$) is reduced.
\end{proof}

\begin{definition}
    A \textit{normal superring} is a Noetherian superring satisfying $(R_1)$ and $(S_2)$.
\end{definition}

\begin{lemma}\label{UnderlyingRingNormal}
    If $A$ is normal, then so is $\ol{A}$.
\end{lemma}

\begin{proof}
    Proven in more generality in Lemma \ref{UnderlyingVarietyNormal}.
\end{proof}

\section{Algebraic Supergeometry}
\subsection{Superschemes and supervarieties}
\begin{definition}
    A \textit{superspace} $X = (|X|, \O_X)$ is a pair consisting of a topological space $|X|$ and a sheaf $\O_X$ of superrings such that each stalk $\O_{X,x}$ is a local superring.
\end{definition}

A \textit{morphism of superspaces} is a pair $\phi = (|\phi|, \phi^\#) : X \to Y$ for a continuous map $|\phi| : |X| \to |Y|$ and a morphism of $\phi^\# : \O_Y \to |\phi|_* \O_X$ of sheaves on $|Y|$ which becomes a local morphism of local superrings on each stalk.

We say that the superspace $X$ is a \textit{superscheme} if $(|X|, (\O_X)_\0)$ is a scheme and the odd part $(\O_X)_\1$ of the structure sheaf is a quasi-coherent sheaf of $(\O_X)_\0$-modules. An \textit{open subscheme} of $X$ is given by the restriction of $\O_X$ to an open set $U \subseteq |X|$, and a \textit{closed subscheme} is an equivalence class of closed immersions in the usual sense.

Superschemes are locally modeled on \textit{affine superschemes} $\Spec A = (|\Spec A_\0|, \tilde A)$ and form a full subcategory of superspaces. Supervarieties will likewise form a full subcategory of superschemes, but we first need the following definitions.

\begin{definition}
    Let $X$ be a superscheme.
    \begin{enumerate}
        \item $X$ is \textit{irreducible} if $|X|$ is irreducible as a topological space.
        \item $X$ is \textit{reduced} if $\O_X(U)$ is a reduced superring for each affine open $U$.
        \item $X$ is \textit{integral} if $\O_X(U)$ is an integral superdomain for each affine open $U$.
        \item $X$ is \textit{GFR} if there is a collection of points $\{\eta\}$ whose closure is $X$ and for which each $\O_{X,\eta}$ is FR.
        \item $X$ is \textit{GFRR} if it is GFR and reduced.
    \end{enumerate}
\end{definition}

\begin{lemma}
Let $X$ be a superscheme.
\begin{enumerate}
    \item $X$ is integral if and only if it is irreducible and reduced.
    \item $X$ is reduced if and only if $\O_{X,x}$ is a reduced superring for all points $x \in |X|$.
\end{enumerate}
\end{lemma}

\begin{proof}
    \begin{enumerate}
        \item Follows from the affine case in Lemma \ref{SuperfieldIntegralReducedLemma}(b) and by passing to the reduced induced closed subschemes for the purpose of using e.g.\ \cite[\href{https://stacks.math.columbia.edu/tag/01OM}{Tag 01OM}]{stacks-project}.
        \item Immediate from Lemma \ref{SuperfieldIntegralReducedLemma}(c). \qedhere
    \end{enumerate}
\end{proof}

The importance of considering only those opens $U$ which are affine in the definition of reduced and integral superschemes is highlighted in the following example.

\begin{example}\label{WeirdP1Example}
    Consider the non-affine superscheme $X$ obtained from gluing the affines $U_1 = \Spec \CC[x, x\xi_1, \xi_2, \xi_1\xi_2]$ and $U_2 = \Spec \CC[x^{-1}, \xi_1, \xi_2]$ along their common localization $\Spec \CC[x^{\pm 1}, \xi_1, \xi_2]$. Then the global coordinate superring $\O_X(X) \cong \CC[\xi_2, \xi_1\xi_2]$ is not reduced, but $X$ certainly is. Geometrically, $U_1$ can be realized as a certain orbit closure of the natural action of $GL(1|2)$ on $S^2 (\CC^{1|2})^*$, appearing in \cite{Sherman}.

    This supervariety is depicted in Figure \ref{fig:fuzzyCircle} as having a non-reduced thick point at $0 \in U_1$. Such behavior is common in the world of supervarieties, when two odd directions interact to produce this type of singularity. We propose that such a singularity be called a ``pill," in analogy with the small balls of lint that occur on surfaces of fabric. In both the mathematical and textile contexts, pills occur when fuzz interacts to produce a small bump.

    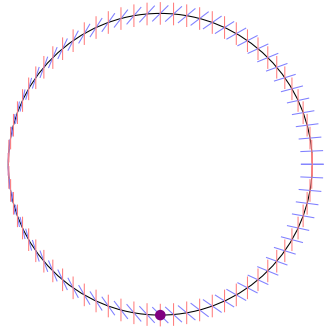
\begin{figure}
        \centering
        \begin{tikzpicture}
            \draw[black] (0,0) circle [radius=2];
            
            \foreach \angle in {0,5,...,360} {
                \draw[blue!50] ({2*cos(\angle)},{2*sin(\angle)}) -- ++({0.15*cos(\angle/2)},{0.15*sin(\angle/2)});
                \draw[blue!50] ({2*cos(\angle)},{2*sin(\angle)}) -- ++({-0.15*cos(\angle/2)},{-0.15*sin(\angle/2)});
                \draw[red!50] ({2*cos(\angle)},{2*sin(\angle)}) -- ++(0,0.15);
                \draw[red!50] ({2*cos(\angle)},{2*sin(\angle)}) -- ++(0,-0.15);
            }
    
            \fill[violet] (0,-2) circle (0.07);
        \end{tikzpicture}
        \caption{A real picture of $X$, with blue fuzz for the $\xi_1$ direction, red fuzz for the $\xi_2$ direction, and a purple ``thick point" at $0 \in U_1$.}
        \label{fig:fuzzyCircle}
    \end{figure}
    
\end{example}

The meanings of separated and finite type are the same as in the usual setting.

\begin{definition}
An \textit{algebraic supervariety} is a GFRR separated superscheme of finite type over an algebraically closed field.
\end{definition}

This definition is designed to mimic the commonly-used definition of a variety found in e.g.\ \cite{Hartshorne}, although we allow for reducible supervarieties by saying ``reduced" instead of ``integral". If we drop the ``GFRR" condition and add ``integral", it becomes equivalent to the formulation given in \cite{Sherman}, although our notion of integrality differs slightly. Here, we consider only those supervarieties which are generically fermionically regular in order to guarantee the existence of smooth points.

As usual, a closed subvariety is a closed subscheme which is itself a supervariety, and an open subvariety is an open subscheme of $X$, viewed as a supervariety. A morphism of supervarieties is a morphism of superschemes between supervarieties.

\begin{lemma}\label{RestrictionMapsInjectiveLemma}
    Let $U \subseteq V \subseteq |X|$ be open subsets of an irreducible supervariety $X$. Then the restriction map $\O_X(V) \to \O_X(U)$ is injective.
\end{lemma}
\begin{proof}
    Let $\eta$ be the generic point. The map $\O_X(V) \to \O_{X,\eta}$ is injective, and factors through $\O_X(V) \to \O_X(U)$.
\end{proof}

\begin{proposition}
    The reduced induced subscheme structure on $X$ is given locally by the underlying rings $\ol{A} = A/J_A$. 
\end{proposition}

\begin{proof}
    It suffices to show that $J_A = \{a \in A \mid a_\p \in \p A_\p \text{ for all } \p \in \Spec A\}$, where $a_\p$ denotes the image of $a$ in $A_\p$. Note that $\supseteq$ holds by definition of $J_A$ using the minimal prime ideal $\p = J_A$. On the other hand, if $a \in J_A$, then $a$ is nilpotent and so $a_\p$ belongs to every prime ideal of $A_\p$.
\end{proof}

We will refer to this closed subvariety of $X$ as the \textit{underlying variety} $X_\0$, and the defining ideal as $\J_X$. Using \cite{Zubkov} and \cite{Hartshorne} Exercise III.3.1, we see that $X$ is affine if and only if $X_\0$ is affine.

\begin{example}
    Let $X$ be as in Example \ref{WeirdP1Example}, so $X_\0 \cong \PP^1$. If we instead consider the closed subscheme of $X$ corresponding to the smaller ideal sheaf $((\O_{X})_1) \subsetneq \J_X$, then we obtain $\PP^1$ with a thick point at $0 \in U_1$, a nonreduced scheme depicted in Figure \ref{fig:fuzzyCircleUnderlyingVariety}(A). It is therefore critical to use the larger ideal sheaf $\J_X$ if we want to genuinely obtain an underlying variety from a supervariety.

    \begin{figure}
        \centering
        
        \begin{subfigure}{.5\textwidth}
          \centering
          \begin{tikzpicture}
            \draw[black] (0,0) circle [radius=2];
            \fill[violet] (0,-2) circle (0.07);
          \end{tikzpicture}
          \caption{Using the ideal sheaf $((\O_X)_\1)$}
        \end{subfigure}%
        \begin{subfigure}{.5\textwidth}
          \centering
          \begin{tikzpicture}
            \draw[black] (0,0) circle [radius=2];
          \end{tikzpicture}
          \caption{The underlying variety}
        \end{subfigure}
        
        \caption{Two different closed subschemes of $X$ as in Example \ref{WeirdP1Example}}
        \label{fig:fuzzyCircleUnderlyingVariety}
    \end{figure}
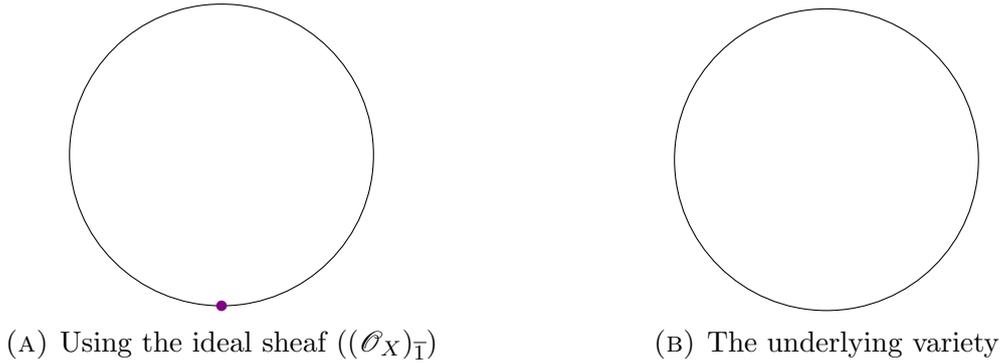
\end{example}

\subsection{Normal supervarieties}

\begin{definition}
    Let $X$ be a Noetherian superscheme. We say $X$ is $(R_k), (S_k)$, or normal if its local rings $\O_{X,x}$ all admit the respective property.
\end{definition}

\begin{definition}
    Let $X$ be a Noetherian superscheme and $Z$ a closed subset of $|X|$ of codimension at least 2. We say $X$ satisfies the \textit{Hartogs property} at $Z$ if the restriction map $\O_X(X) \to \O_X(X \backslash Z)$ is an isomorphism. Further, we say that $X$ satisfies the Hartogs property (irrespective of any closed subset) if it satisfies the Hartogs property at all closed subsets of codimension at least 2.
\end{definition}

We record the following mild generalization of a standard lemma.

\begin{lemma}\label{HartogsLemma}
    Let $X$ be a Noetherian superscheme. Then $X$ satisfies $(S_2)$ if and only if it satisfies $(S_1)$ and the Hartogs property.
\end{lemma}

\begin{proof}
    Since all relevant properties are local, we may assume $X = \Spec A$ is affine. Note that $\O_X$ is a coherent sheaf of $(\O_X)_\0$-modules, and $(X, (\O_X)_\0)$ is an affine scheme, so $H^1(X,\O_X) = 0$. Now let $\p \in |\Spec A|$ be a point and $Z$ its closure. If $\height \p \leq 1$ then we are finished because $(S_2)$ implies $(S_1)$. Now suppose $\height \p \geq 2$. The long exact sequence in cohomology with supports truncates as
    \begin{center}
    \begin{tikzcd}[column sep = small]
        0 \arrow[r] & H^0_Z(X,\O_X) \arrow[r] & H^0(X,\O_X) \arrow[r] & H^0(X\backslash Z, \O_X|_{X\backslash Z}) \arrow[r] & H^1_Z(X,\O_X) \arrow[r] & 0
    \end{tikzcd}
    \end{center}
    so $X$ satisfies the Hartogs property at $Z$ if and only if $H^i_Z(X,\O_X) = 0$ for $i=0,1$. Using \cite{Hartshorne} Exercises III.3.3-4, this condition is equivalent to $\Ext^i(A_\p/\p, A_\p) = 0$ for $i=0,1$, which is $(S_2)$.
\end{proof}

\begin{lemma}\label{UnderlyingVarietyNormal}
    If $X$ is a normal supervariety, then so is $X_\0$.
\end{lemma}
\begin{proof}
    The condition $(R_1)$ is immediate, and implies $(S_1)$. For $(S_2)$, we note that the Hartogs property of Lemma \ref{HartogsLemma} implies the Hartogs property for $X_\0$, since $\ol{\O_X(X)} \to \ol{\O_X(X\backslash Z)}$ is also an isomorphism. The claim then follows by Lemma $\ref{HartogsLemma}$.
\end{proof}

\subsection{Algebraic supergroups and their actions}
The material in this section is well-known and provided for reference.

\begin{definition}
    An \textit{affine algebraic supergroup} is any of the following equivalent notions:
    \begin{enumerate}
        \item The spectrum of a finitely-generated Hopf superalgebra $A$ which is an integral superdomain
        
        \item A group object $G$ in the category of affine supervarieties
        
        \item An affine supervariety $G = \Spec A$ whose functor of points $G(B) := \Hom(A,B)$ is group-valued
        
        \item A pair $G = (G_\0, \g)$ where $G_\0$ is an ordinary affine algebraic group, $\g = \g_\0 \oplus \g_\1$ is a finite-dimensional Lie superalgebra such that $\Lie G_\0 = \g_\0$, and $\g$ is equipped with a $G_\0$-module structure such that the corresponding $\g_0$-module is adjoint (called a \textit{Harish-Chandra pair})
    \end{enumerate}
\end{definition}

We use these notions (which can be found in \cite{CCF}, Sections 7 and 11) interchangeably throughout this paper. The requirement that the Hopf superalgebra $R$ of part (a) be integral is needed only to guarantee irreducibility of the supergroup. If we relax the assumption that a supervariety must be irreducible (and accordingly adjust for the possibility of multiple generic points), then this requirement may be dropped, and definitions (a)-(d) remain equivalent.

The notation in the definition of a Harish-Chandra pair is intentional, as $G_\0$ is the underlying algebraic group of $G$ and $\g$ is the Lie superalgebra (i.e.\ tangent space at the identity) of $G$.

\begin{definition}
    A (left) action of a supergroup $G$ on a supervariety $X$ is a morphism $a : G \times X \to X$ of supervarieties satisfying all the usual properties.
\end{definition}

If $y \in X(\kk)$ is a closed point, then there is an orbit map $a_y : G \to X$ given by $a \circ (id_G \times i_y)$, where $i_y : y \to X$ is the natural inclusion. The \textit{stabilizer subgroup} $\Stab_G(y)$ is a closed subgroup (see \cite{Fioresi}) of $G$ defined as the fiber $a_y^{-1}(x)$. The orbit map $a_y : G \to X$ factors through the natural surjection $G \to G/\Stab_G(y)$, and $G/\Stab_G(y) \to X$ is an immersion of supervarieties. If this immersion is open, we refer to it as an open orbit of $G$.

The orbit closure of $y$, denoted $\cl(G \cdot y)$, is the closed subscheme of $X$ determined by the sheaf of ideals given by the kernel of the pullback map $\O_X \to (a_y)_* \O_G$. It is seen in \cite{Sherman} that $\cl(G \cdot y)$ is a supervariety with open orbit isomorphic to $G / \Stab_G(y)$.

\section{Toric supervarieties}

\subsection{Algebraic supertori}

\begin{definition}
    An \textit{algebraic supertorus} is an algebraic supergroup $T = (T_\0, \t)$ such that $T_\0 \cong (\kk^\times)^n$ is an ordinary algebraic torus lying in the center of $T$.
\end{definition}

This definition is motivated by the theory of Cartan subgroups of quasireductive supergroups as in \cite{SerganovaQred}. A \textit{quasireductive} supergroup $G$ is one whose even part $G_\0$ is reductive. In this setting, one chooses a Cartan subgroup (i.e.\ maximal torus) $T_\0$ of $G_\0$, and defines $T$ as the centralizer of $T_\0$ in $G$. Thus, an algebraic supertorus is an algebraic supergroup which occurs as a Cartan subgroup of a quasireductive supergroup. The representation theory of supertori has been studied in \cite{GSS} and \cite{Shibata}.

From the definition of a Harish-Chandra pair, one sees that a supertorus is determined by an ordinary algebraic torus $T_\0$ and a super vector space $\t = \t_\0 \oplus \t_\1$ (such that $\t_0 = \Lie T_\0$) equipped with a linear map $S^2 \t_\1 \to \t_\0$.

We write $M$ and $N$ for the weight and coweight lattices, respectively, of $T_\0$. Conversely, if $N$ is a lattice, we write $T_N$ for the algebraic torus whose coweight lattice is $N$. The perfect pairing between $M$ and $N$ will be written as 
$$\langle -, - \rangle : M \times N \to \ZZ,$$
and we note that it extends $\CC$-linearly to $M \times \t_\0$ since $\t_\0 \cong N_\CC := N \otimes_\ZZ \CC$.

The algebra of regular functions on $T_\0$ is isomorphic to the group algebra $\CC[M]$, and likewise $\CC[T] \cong \CC[M] \otimes_\CC \bigwedge \t_\1^*$. In coordinates, we will often fix compatible isomorphisms
\begin{align*}
    T_N &\cong (\kk^\times)^{p} \\
    \kk[T_N] &\cong \kk[t_1^{\pm 1}, ..., t_{p}^{\pm 1}], \\
    \kk[T] &\cong \kk[t_1^{\pm 1}, ..., t_{p}^{\pm 1}, \xi_1, ..., \xi_{q}], \\
    \t &\cong \kk\{ x_1, ..., x_{p} \} \oplus \kk\{ \theta_1, ..., \theta_{q} \}.
\end{align*}
where we view the $x_i$ (respectively, $\theta_i$) $\in \Der_{\CC[T]}(\CC[T], \CC)$ as even (respectively, odd) point derivations at the identity such that $x_i(t_j) = \delta_{ij} = \theta_i(\xi_j)$.

We then choose $x_{ij} = (x_{ij})_1 x_1 + ... + (x_{ij})_{p} x_{p} \in \t_\0$ for $i,j = 1, ..., q$ such that $x_{ij}+x_{ji} = [\theta_i,\theta_j]$ for all such $i$ and $j$. In some cases we will use $x_{ij} = \frac{1}{2}[\theta_i,\theta_j]$, although in some situations it will be more convenient to make different choices. We will express the Hopf superalgebra structure on $\CC[T]$, the group law of $T$, and the action of $\t$ on $\CC[T]$ in terms of these choices, but we stress that any two choices are equivalent via automorphisms of $T$.

The Hopf superalgebra structure on $\CC[T]$ is given by the rules
\begin{align*}
    \Delta t^m &= t^m \otimes t^m \left(1 - \sum_{i,j=1}^{q} \langle m, x_{ij} \rangle \xi_i \otimes \xi_j \right) & S t^m &= t^{-m} & \epsilon t^m &= 1 \\
    \Delta \xi_i &= \xi_i \otimes 1 + 1 \otimes \xi_i &
    S \xi_i &= -\xi_i & \epsilon \xi_i &= 0
\end{align*}
where $t^m \in \CC[M]$ is the function corresponding to the character $m \in M$, and the group law on $A$-points is \begin{align*}
    &(a_1, ..., a_{p} \mid \alpha_1, ..., \alpha_{q}) \cdot (b_1, ..., b_{p} \mid \beta_1, ..., \beta_{q}) \\
    &= \left(a_1 b_1 \left(1- \sum_{i,j=1}^{q} (x_{ij})_1 \alpha_i \beta_j \right), ..., a_{p} b_{p} \left(1- \sum_{i,j=1}^{q} (x_{ij})_{p} \alpha_i \beta_j\right) \; \Bigg| \; \alpha_1 + \beta_1, ..., \alpha_{q} + \beta_{q} \right)
\end{align*}
for $a_i, b_i \in (A_\0)^\times$ and $\alpha_i, \beta_i \in A_\1$, and 

The left action of $T$ on itself induces the left regular representation of $T$ on $\CC[T]$. This action differentiates to one of $\t$ on $\CC[T]$, wherein
\begin{align}
\begin{split}
    x_i &\mapsto -t_i \d{t_i} \\
    \theta_i &\mapsto \sum_{j=1}^{q} \xi_j \left\langle t \d{t},  x_{ij}  \right\rangle - \d{\xi_i}
\end{split}\label{eqn:LeftDerivations}
\end{align}
and where we have written $\left\langle t \d{t},  x_{ij}  \right\rangle$ to signify $(x_{ij})_1 t_1 \d{t_1} + ... + (x_{ij})_{p} t_{p} \d{t_{p}}$, by abuse of the notation $\langle -,-\rangle : M \times N \to \ZZ$. For $x \in \t$, we write $x$ for the corresponding derivation. In general when we refer to ``the" action of $T$ (or $\t$) on $\CC[T]$, we will mean the left regular representation.  However, we will also later need the right regular representation, for which we will write $x^{\r}$, so we record this action of $\t$ here:
\begin{align*}
\begin{split}
    x_i^\r &\mapsto t_i \d{t_i} \\
    \theta_i^\r &\mapsto \sum_{j=1}^{q} \xi_j \left\langle t \d{t},  x_{ji}  \right\rangle + \d{\xi_i}
\end{split}
\end{align*}

When $H$ is a subgroup of $T$, we will write $\CC[T]^H$ (or $A^H$) for the $H$-invariants of $\CC[T]$ (or of $A \subseteq \CC[T]$) under the right regular representation, so that $\CC[T]^H \cong \CC[T/H]$. Likewise, if $x$ is a vector field on $A$, we write $A^x$ for its kernel.

\begin{example}\label{Q(1)nExample}
    A recurring example in this paper will be the supertorus $Q(1)^n$, for which $p=q=n$ and $\frac{1}{2}[\theta_i,\theta_j] = \delta_{ij} x_i$.
\end{example}

\subsection{Introduction to toric supervarieties}

We now begin our treatment of toric supervarieties. While the following definition requires that the supertorus action be a left action, all the same theory (with a few sign changes, noted above) applies if we switch to right actions. Let us first recall the definition of a toric variety.

\begin{definition}
    A \textit{toric variety} is an irreducible normal variety $X$ equipped with a left action by an algebraic torus $T$ and a $T$-equivariant open immersion $T \to X$.
\end{definition}

Notice that if $X$ is any (irreducible, normal) variety with an open orbit of a torus $T$, then the open orbit $T/H$ itself is an algebraic torus and $X$ is a toric variety for $T/H$. In this sense, any such variety with an open orbit of a torus can be modeled as a toric variety according to the above definition.

If we wish to generalize this notion to the supergeometric setting, we must therefore contend with the realm of homogeneous spaces $T/H$ for a supertorus $T$. Since $T_\0$ is central in $T$, any purely even subgroup $H_\0 \subseteq T$ is a normal subgroup, and so $T/H_\0$ is a supertorus. The same is not true for an arbitrary supergroup $H$; in fact, $T/H$ need not be a supergroup at all!

\begin{example}
    Consider the supergroup $T = Q(1)^2 := Q(1) \times Q(1)$, with subgroup $H = Q(1)$ embedded diagonally. It is routine to check that $H$ is not a normal subgroup, and so $T/H$ does not inherit the structure of a supergroup.
\end{example}

Consequently, in our definition, we make the compromise that the open orbit of a toric supervariety is isomorphic merely to a homogeneous space $T/H$ for a purely odd (and therefore abelian) subgroup $H$.

\begin{definition}
    A \textit{toric supervariety} is an irreducible normal supervariety $X$ equipped with a left action by a supertorus $T$, and a $T$-equivariant map $T \to X$ which factors through an open immersion of a homogeneous space $T/H$. We will assume $H_\0=1$. If $H=1$, then the toric supervariety will be called \textit{faithful}.
\end{definition}

Occasionally we will wish to drop the normality criterion, at which times it will be clarified that the supervariety is not necessarily normal.

\begin{lemma}
    Let $X$ be a (not necessarily normal) toric supervariety with torus $T$. Then the underlying variety $X_\0$ is a (not necessarily normal) toric variety with torus $T_\0$.
\end{lemma}

\begin{proof}
    The action of $T_\0$ on $X$ preserves the defining ideal $\J_X$ of $X_\0$.
\end{proof}

\subsection{Structure of affine toric supervarieties}
The problem of understanding toric supervarieties reduces to the affine case, so we will first attempt to characterize which superalgebras $A$ can occur as their regular functions.

\begin{lemma}\label{ATVCharacterizationLemma}
    Let $X = \Spec A$ be an affine supervariety. Then:
    \begin{enumerate}
        \item $X$ is toric (but not necessarily normal) for a quotient of the supertorus $T$ if and only if $A$ is isomorphic to a finitely-generated subalgebra of $\CC[T]$ which is also a $\t$-subrepresentation.
        \item If $X$ is toric (but not necessarily normal) for the supertorus $T$, then it is faithful if and only if there is a non-nilpotent $T_\0$-weight vector $f \in A$ such that $A[f^{-1}] \cong \CC[T]$.
    \end{enumerate}
\end{lemma}
\begin{proof}
    First assume $A \subseteq \CC[T]$ is such a subalgebra as in (a). The inclusion $A \to \CC[T]$ dualizes to map $T \to X$ which factors through an open immersion $U \to X$ with $T \to U$ surjective. Because $A$ is a subrepresentation of $\CC[T]$, all the above maps of supervarieties are $T$-equivariant. Hence $U \cong T/H$ for some subgroup $H \subseteq T$, so $X$ is a toric supervariety for $T/H_\0$. If there is a suitable weight vector $f \in A$ as in (b), then $U \cong T$ is the basic open $D(f)$ and so $X$ is faithful.

    Conversely, suppose $X$ is toric with open orbit isomorphic to $T/H$. By Lemma \ref{RestrictionMapsInjectiveLemma}, $A \to \CC[T/H] = \CC[T]^H \to \CC[T]$ is injective. Moreover, the action of $T$ on $X$ becomes an action of $T$ on $\CC[X]$, so the conditions of (a) are satisfied. Now if $X$ is faithful, consider the underlying variety $X_\0 \supseteq T_\0$ with reduced coordinate algebra $\ol{A} \subseteq \CC[T_\0]$. Here, there is a $T_\0$-weight vector $t^m \in \ol{A}$ such that $\ol{A}[t^{-m}] = \CC[T_\0]$; i.e.\ $T_\0 = D(t^m)$. Let $f \in A$ lie in the preimage of $x^m$ under the projection $A \to \ol{A}$. Then the open subvariety $D(f)$ has underlying variety $T_\0$ and so $A[f^{-1}] \cong \CC[T]$.
\end{proof}

Henceforth consider only those toric supervarieties which are normal. By Lemma \ref{UnderlyingVarietyNormal}, this means the underlying toric varieties are normal, so the affines are determined by strongly convex, rational, polyhedral cones $\sigma \subset N_\RR := N \otimes_\ZZ \RR$. When we say ``cone" we implicitly assume all the prior adjectives.

More generally, a toric variety $X_\Sigma$ for the torus $T_N$ induces a fan $\Sigma$ in $N_\RR$, i.e.\ a collection of cones which contains all the faces of each of its cones, and such that the intersection of any two cones in $\Sigma$ is a face of each. If $\tau$ is a face of $\sigma$, we write $\tau \leq \sigma$, so that $\kk[S_\tau] \supseteq \kk[S_\sigma]$ at the level of semigroup algebras, where $S_\sigma = \sigma^\vee \cap M$ is the semigroup of lattice points in the dual cone $\sigma^\vee$. For $m, m' \in S_\sigma$, we write $m' \leq_\sigma m$ if $m-m' \in S_\sigma$. We refer to \cite{CLS} for any missing exposition on the subject.

We will use the symbol $\bullet$ to refer to the 0-dimensional cone, and $\rho$ to refer to a ray. The set of rays belonging to the fan $\Sigma$ (respectively, the cone $\sigma$) will be denoted by $\Sigma(1)$ (respectively, $\sigma(1)$). For an arbitrary cone $\sigma \in \Sigma$, we write $A_\sigma$ for the coordinate superalgebra of an affine toric supervariety such that $\ol{A_\sigma} = \CC[S_\sigma]$.

If $\Spec A$ is an affine toric supervariety where $A=A_\sigma$ has underlying toric variety $\Spec \kk[S_\sigma]$, we write $A_\tau$ for the coordinate superalgebra of the affine toric open subvariety corresponding to the face $\tau \leq \sigma$. Moreover, if $m \in S_\sigma$, we write $L_{A_\sigma}(m)$ for the $-m$-weight space of $A_\sigma$ as a $T_\0$-module. (The apparent sign discrepancy occurs because $T_\0$ acts by the character $-m$ on the monomial $t^m$.) Sometimes we will suppress the subscript and write $L(m)$ if the meaning is clear from context.

The notation $L(m)$ is intended to be evocative of, but not confused with, the typical notation for the irreducible weight modules of a torus. However, we stress that $L_{\CC[T]}(m)$ is not usually irreducible for a supertorus $T$.

\begin{lemma}\label{RaySuperalgebraGeneratorsLemma}
    Let $A_\rho$ be the coordinate superalgebra of a (not necessarily normal) affine toric supervariety corresponding to a ray $\rho$. Then $A_\rho$ may be written as $$\kk \left[ t^{m_1}(1+...), \left(t^{m_2}(1+...) \right)^{\pm 1}, ..., \left(t^{m_{p}}(1+...) \right)^{\pm 1}, t^{\ell_1 m_1}(...), ..., t^{\ell_r m_1}(...) \right]$$ for $m_i$ a basis of $M$ and $\ell_i \in \ZZ_{\geq 0}$, where any ellipsis within parentheses, $(1+...)$ or $(...)$, refers to nilpotent terms in $L(0)$.
\end{lemma}
\begin{proof}
    We know $\ol{A_\rho} \cong \CC[t^{m_1}, t^{\pm m_2}, ..., t^{\pm m_{p}}]$. If $f,g \in A_\rho$ are lifts of inverse elements in $\ol{A_\rho}$, then $fg = 1+\eta$ for some $\eta \in J_{A_\rho}$, so $f(1-\eta-\eta^2-...) = g^{-1}$ and hence we may assume that inverses lift to inverses. The rest follows since we may generate $A_\rho$ by $T_\0$-weight vectors.
\end{proof}

\begin{proposition}\label{NormalATSVProp}
    Let $X = \Spec A$ be a (not necessarily normal) affine toric supervariety. Then $X$ is normal if and only if all of the following conditions hold:
    \begin{enumerate}
        \item $\ol{A} \cong \kk[S_\sigma]$ for some cone $\sigma \subseteq N_\RR$
        \item $A_\sigma = \bigcap_{\rho \in \sigma(1)} A_\rho$
        \item $A_\rho$ is FR for all $\rho$
    \end{enumerate}
\end{proposition}

\begin{proof}
    First assume $X$ is normal, so $X_\0$ is too and hence (a) holds. We also have (b) by the Hartogs property of Lemma \ref{HartogsLemma}, so it remains to prove (c). Using the expression for $A_\rho$ obtained in Lemma \ref{RaySuperalgebraGeneratorsLemma}, consider the height-1 prime ideal generated by nilpotent elements and $t^{m_1}(1+...)$. Localization at this prime preserves the nilpotent generators, so they must form a regular sequence. Then $A_\rho$ is isomorphic to the exterior algebra of the span of this regular sequence over $\ol{A_\rho}$, so (c) follows.

    For the converse, we note that (a) and (c) together imply that $U = \bigcup_{\rho \in \sigma(1)} \Spec A_\rho$ is regular and therefore normal. Then the singular locus of $X$ has codimension at least 2, so $(R_1)$, and hence also $(S_1)$, holds. For $(S_2)$, it suffices by Lemma \ref{HartogsLemma} to show that $X$ admits the Hartogs property. Let $Z$ be a closed subset of $|X|$ with codimension at least 2, so $Z \cap U$ is likewise a closed subset of $U$ with codimension at least 2. Since $U$ is smooth and hence $(S_2)$, we have
    \begin{center}
    \begin{tikzcd}
        \O_X(X) \arrow[r, "\sim"] \arrow[hookrightarrow,d] & \O_X(U) \arrow[d, "\sim"] \\
        \O_X(X\backslash Z) \arrow[hookrightarrow,r] & \O_X(U \backslash Z)
    \end{tikzcd}
    \end{center}
    so indeed $X$ has the Hartogs property at $Z$. 
\end{proof}

Thus, to determine an affine toric supervariety, one needs only a cone $\sigma \subset N_\RR$ and a suitable FR superalgebra $A_\rho$ for each of its rays $\rho \in \sigma(1)$. In the following section, we propose a method of extracting classifying data from the coordinate superalgebras.

\subsection{Decorated fans}

Let $X$ be a toric supervariety (for the supertorus $T$) with open orbit $T/H$, and suppose $X_\0 \cong X_\Sigma$. Using Lemma \ref{ATVCharacterizationLemma}(a), we may assume that each $A_\sigma$ (for $\sigma \in \Sigma$) is a subalgebra of $\CC[T]$ in such a way that $A_\sigma$ localizes to $A_\tau$ whenever $\tau \leq \sigma$. For each cone $\sigma \in \Sigma$ and each character $m \in S_\sigma$, let $A_{\sigma,m}$ be the subalgebra of $A_\sigma$ generated by $L_{A_\sigma}(m')$ for all $m' \leq_\sigma m$. Then $A_{\sigma,m}$ also satisfies the condition of Lemma \ref{ATVCharacterizationLemma}(a), so we may define $H_{\sigma,m}$ as the stabilizer of the image of 1, i.e.\ the subtorus of $T$ for which the $T$-action on $\Spec A_{\sigma,m}$ has open orbit $T/H_{\sigma,m}$.

The following lemma says that the choice of the point 1 was largely irrelevant.

\begin{lemma}
    Let $T$ be a supertorus, $Y$ a homogeneous space for $T$, and $x,y \in Y$ two closed points. Then $\Stab_T x = \Stab_T y$.
\end{lemma}
\begin{proof}
    Let $g \in T_\0$ be such that $g x = y$. Then $\Stab_T y$ can be obtained from $\Stab_T x$ via conjugation by $g$, which is central. Therefore $\Stab_T x = \Stab_T y$.
\end{proof}

We write $\h := \Lie H$ and $\h_{\sigma,m} := \Lie H_{\sigma,m}$. Moreover, for $m \in S_\sigma$, write $\sigma_m = \{x \in \sigma \mid \langle m, x \rangle = 0\}$ for the face of $\sigma$ consisting of its intersection with the hyperplane $\ker m$. We also write $N_\sigma$ for the sublattice of $N$ generated by $\sigma \cap N$.

\begin{lemma}\label{TSVDataLemma}
    Using the above notation, we have:
    \begin{enumerate}
        \item $\ol{A_{\sigma,m}} \cong \CC[\sigma^\vee \cap \sigma_m^\perp \cap M]$

        \item $(H_{\sigma,m})_\0 = T_{N_{\sigma_m}}$

        \item $(\h_{\sigma,m})_\1 = \{\theta \in \t_1 \mid \theta^\r L_{A_\sigma}(m) = 0\}$

        \item $H_{\sigma,m} \supseteq H_{\sigma,m'}$ whenever $m \leq_\sigma m'$

        \item $\h_{\sigma,m} \supseteq \sum_{\rho \in \sigma(1)} \h_{\rho,m}$
    \end{enumerate}
\end{lemma}

\begin{proof}
\begin{enumerate}
    \item This is because $\sigma^\vee \cap \sigma_m^\perp \leq \sigma^\vee$ is the dual face of $\sigma_m \leq \sigma$, and $m$ belongs to the relative interior of $\sigma_m^\perp$.

    \item The open $T_\0$-orbit of $\Spec \ol{A_{\sigma,m}}$ has coordinate algebra $\CC[\sigma_m^\perp \cap M] \cong \CC[(N/N_{\sigma_m})^*]$ by part (a).

    \item Since $\CC[T/H_{\sigma,m}] \cong \CC[T_\0/T_{N_\tau}] \otimes \bigwedge (\t_\1 / (\h_{\sigma,m})_\1)^*$ is a localization of $A_{\sigma,m}$, we have $(\h_{\sigma,m})_\1 = \{\theta \in \t_1 \mid \theta^{\r} L_{\CC[T/H_{\sigma,m}]}(0) = 0\} = \{\theta \in \t_1 \mid \theta^{\r}  L_{A_\sigma}(m) = 0\}$.
    
    \item Immediate from construction of $A_{\sigma,m} \subseteq A_{\sigma,m'}$.
    
    \item We have $L_{A_\sigma}(m) = \bigcap_{\rho \leq \sigma} L_{A_\rho}(m)$ by Proposition \ref{NormalATSVProp}(b), so the statement holds for the odd part by (c). Equality of the even parts follows from (b). \qedhere
\end{enumerate}
\end{proof}

If $\rho \in \Sigma(1)$, we write $u_\rho$ for the primitive generator of the ray $\rho$. That is, $u_\rho$ is the first nonzero lattice point along $\rho$. The following lemma is an immediate consequence of Lemma \ref{TSVDataLemma}(d).

\begin{lemma}
    Let $\rho \in \Sigma(1)$. Then $\h_{\rho,m}$ depends only on the nonnegative integer $\langle m, u_\rho \rangle$.
\end{lemma}

For $i \geq 0$, we write $V_{\rho,i}$ (or $V_{\rho,m}$) for any subspace $(\h_{\rho,m})_\1 \subseteq \t_\1$ such that $\langle m, u_\rho \rangle = i$. The prior lemma asserts that the choice of $m$ is immaterial. We collect some results about the subspaces $V_{\rho,i}$ here:

\begin{lemma}\label{lemma:DecorationCompatibility}
Let $V_{\rho,i}$ be as above. Then:
\begin{enumerate}
    \item $V_{\rho,0} \supseteq V_{\rho,1} \supseteq V_{\rho,2} \supseteq ...$
    \item The chain $V_{\rho,0} \supseteq V_{\rho,1} \supseteq V_{\rho,2} \supseteq ...$ stabilizes at $\h$.
    \item If $\rho, \rho' \in \sigma(1)$ and $m \in \sigma^\vee \cap M$, then $[V_{\rho,m}, V_{\rho',m}] \subseteq \t_{N_{\sigma_m}}$. In particular, $[V_{\rho,0}, V_{\rho,0}] \subseteq \t_{N_\rho}$ and $[V_{\rho,1}, V_{\rho,1}] = 0$.
    \item Either $[V_{\rho,0}, V_{\rho,0}] =0$ or $\codim(\h, V_{\rho,0})=1$.
\end{enumerate}
\end{lemma}

\begin{proof}
\begin{enumerate}
    \item Lemma \ref{TSVDataLemma}(d).
    \item By Lemma \ref{RaySuperalgebraGeneratorsLemma}, $A_\rho \to A_\bullet = \CC[T/H]$ is a localization via the inversion of a single even weight vector $f = t^{m_1}(1+...) \in A_\rho$. Hence $$L_{\CC[T/H]}(0) = \sum_{n \geq 0} L_{A_\rho}(nm_1)[f^{-n}],$$ so there is $n \geq 0$ for which $$L_{\CC[T/H]}(0) = L_{A_\rho}(nm_1)[f^{-n}].$$ Therefore $V_{\rho,n} = \h$.
    \item Lemma \ref{TSVDataLemma} parts (b) and (e).
    \item Suppose $\codim(\h,V_{\rho,0}) \geq 2$ and $[V_{\rho,0}, V_{\rho,0}] \neq 0$, so we may write $$A_\rho = \CC[t_1(1+...), (t_2(1+...))^{\pm 1}, ..., (t_p(1+...))^{\pm 1}, \xi_1, ..., \xi_r, t_1^{\ell_{r+1}}(\xi_{r+1} + ...), ..., t_1^{\ell_q}(\xi_q+...)]$$ and pick $\theta \in V_{\rho,0}$ such that $[\theta,\theta] \neq 0$. We may assume $\theta = \xi_q t_1 \d{t_1} - \d{\xi_q}$. Then for $i=r+1, ..., q-1$, we have $\theta \cdot (t_1^{\ell_i}(\xi_i+...)) = t_1^{\ell_i}(\xi_q \xi_i + \xi_q(...) - \d{\xi_q} (...))$. But $\d{\xi_q} (...)$ has no terms that can cancel with $\xi_q \xi_i$, contradicting that $\xi_q \xi_i$ first appears in $L(\ell_i + \ell_q)$. \qedhere
\end{enumerate}
\end{proof}

In other words, a toric supervariety induces a decorated fan in the following sense:
\begin{definition}
    A \textit{decorated fan} is a tuple $(N,\Sigma,\t,\h,\{V_{\rho,i}\}_{\rho \in \Sigma(1), i \geq 0})$ consisting of a finite-dimensional lattice $N$, a fan $\Sigma$ in $N_\RR$, a quasiabelian Lie superalgebra $\t$ with $\t_\0 = \t_N$, an odd abelian subalgebra $\h \subseteq \t$, and a collection of descending chains $V_{\rho,0} \supseteq V_{\rho,1} \supseteq ...$ of subspaces of $\t_\1$, stabilizing at $\h$.

    Beginning in Definition \ref{def:LETDecoratedFan}, we will additionally impose compatibility conditions on the $V_{\rho,i}$.
\end{definition}

Our objective now is threefold: to determine which decorated fans arise from toric supervarieties, to determine to what extent we can recover a toric supervariety from its decorated fan, and to extract geometric information about a toric supervariety from its decorated fan. A first step towards the latter goal occurs in the following proposition, wherein we write $V_{\sigma,0} = \sum_{\rho \in \sigma(1)} V_{\rho,0}$ for the odd part of the Lie superalgebra $\h_{\sigma,0}$.

\begin{proposition}\label{prop:orbit-stabilizer}
    Let $X$ be a toric supervariety with decorated fan $(N,\Sigma,\t,\h,\{V_{\rho,i}\})$, and let $y \in X$ be a closed point in the $T_\0$-orbit corresponding to $\sigma \in \Sigma$. Then $\Stab_T y = (T_{N_\sigma}, \t_{N_\sigma} \oplus V_{\sigma,0})$.
\end{proposition}

\begin{proof}
    The even part of the orbit is determined entirely classically. Meanwhile, the odd part is $V_{\sigma,0}$ because any odd coordinates $t^m (...)$ in $A_\sigma$ will be fixed at 0 in this orbit unless $m \in \sigma^\perp$.
\end{proof}

\subsection{Toric supervarieties with the same decorated fan}
In many important situations, a decorated fan will correspond to a unique toric supervariety. However, for supertori with large odd abelian subalgebras, this is far from the case.

\begin{example}\label{ex:WildTSV}
    Let $T$ be the abelian supertorus of dimension $(1|4)$, so that its Lie superalgebra acts by the derivations $t\d{t}$ and $\d{\xi_i}$ for $i=1, ..., 4$ on $\CC[T] = \CC[t^{\pm 1}, \xi_1, \xi_2, \xi_3, \xi_4]$. We consider the following two superalgebras:
    \begin{align*}
        A &= \CC[t, t\xi_1, t\xi_2, t\xi_3, t^2\xi_4] \\
        A' &= \CC[t, t\xi_1, t\xi_2, t\xi_3, t^2 (\xi_4+\xi_1\xi_2\xi_3)]
    \end{align*}
    Both $\Spec A$ and $\Spec A'$ admit the same decorations $\t_\1 \supseteq \CC\theta_4 \supseteq 0 \supseteq ...$, and both are isomorphic to $\AA^{1|4}$ as a supervariety. However, the actions of $T$ are decidedly different.
\end{example}

In fact, for particularly large abelian supertori, classification is a wild problem.

\begin{proposition}
    A classification of toric supervarieties contains a classification of triples of skew-symmetric bilinear forms up to simultaneous congruence.
\end{proposition}
\begin{proof}
    Let $T$ be the abelian supertorus of dimension $(3|n+4)$. Write $$\CC[T] = \CC[t_A^{\pm 1}, t_B^{\pm 1}, t_C^{\pm 1}, \xi_1, ..., \xi_n, \xi_A, \xi_B, \xi_C, \xi'],$$ and assume the odd part $\t_\1$ of the Lie superalgebra $\t$ acts by $\theta_i = -\frac{d}{d\xi_i}.$
    
    Let $\sigma$ be the positive orthant, corresponding to the toric variety $\Spec \CC[t_A, t_B, t_C]$. Write $\rho_A, \rho_B, \rho_C$ for its rays, and let 
    \begin{align*}
    V_{\rho_\ell,0} &= \CC\{\theta_1, ..., \theta_n, \theta_\ell, \theta'\} \\
    V_{\rho_\ell,1} &= \CC\{\theta_\ell, \theta'\} \\
    V_{\rho_\ell,2} &= \CC\{\theta'\} \\
    V_{\rho_\ell,3} &= 0
    \end{align*}
    for $\ell=A,B,C$. 
    
    We will attempt to classify toric supervarieties for $T$ inducing the decorated fan $(\ZZ^3, \sigma, \t, 0, \{V_{\rho,i}\})$ specified above. For each $\ell$, suppose $\theta_i\theta_j\theta_k$ kills $L_{A_{\rho_\ell}}(2)$ for any $i,j,k \in \{1, ..., n\}$. Then the possible coordinate superalgebras of $U_{\rho_\ell}$ are
    \begin{align*}
        A_{\rho_A} &= \CC \left[t_A, t_B^{\pm 1}, t_C^{\pm 1}, t_A \xi_1, ..., t_A \xi_n, t_A^2 \xi_A, \xi_B, \xi_C, t_A^3 \left(\xi' + \xi_A \left(\sum_{1 \leq i < j \leq n} a_{ij} \xi_i \xi_j \right) \right) \right] \\
        A_{\rho_B} &= \CC \left[t_A^{\pm 1}, t_B, t_C^{\pm 1}, t_B \xi_1, ..., t_B \xi_n, \xi_A, t_B^2 \xi_B, \xi_C, t_B^3 \left(\xi' + \xi_B \left(\sum_{1 \leq i < j \leq n} b_{ij} \xi_i \xi_j \right) \right) \right] \\
        A_{\rho_C} &= \CC \left[t_A^{\pm 1}, t_B^{\pm 1}, t_C, t_C \xi_1, ..., t_C \xi_n, \xi_A, \xi_B, t_C^2\xi_C, t_C^3 \left(\xi' + \xi_C \left(\sum_{1 \leq i < j \leq n} c_{ij} \xi_i \xi_j \right) \right) \right]
    \end{align*}
    and hence we may write $(a_{ij})$, $(b_{ij})$, $(c_{ij})$ as skew-symmetric matrices. But these matrices depend on a choice of basis of $\CC\{\xi_1, ..., \xi_n\}$, and such changes of basis are equivalent to simultaneous congruence transformations. 
\end{proof}

In summary, one decorated fan can describe many different toric supervarieties. We introduce the following geometric criterion to ensure that the correspondence becomes bijective.

\subsection{Homological regularity}

\begin{propdef}\label{propdef:HR1}
Let $X$ be a toric supervariety with decorated fan $(N,\Sigma,\t,\h,\{V_{\rho,i}\})$. We say $X$ is \textit{homologically regular in codimension 1}, or $(HR_1)$, if it satisfies any of the following equivalent conditions:
\begin{enumerate}
    \item For every $\rho \in \Sigma(1)$, there is a parameterization as in (\ref{eqn:LeftDerivations}) of the derivations by which $\t$ acts on $\CC[T]$ such that
    $$A_\rho = \CC[t_1, t_2^{\pm 1}, ..., t_p^{\pm 1}, \xi_1, ..., \xi_r, t_1^{\ell_{r+1}} \xi_{r+1}, ..., t_1^{\ell_q} \xi_q]$$
    for some positive integers $\ell_i$.
    
    \item For every $\rho \in \Sigma(1)$, $A_\rho$ is the subalgebra of $\CC[T]$ generated by $L_{\CC[T]^W}(n+\rho^\perp)$ for all $n \geq 0$ and all codimension-at-most-1 subspaces $W \subseteq V_{\rho,0}$ containing $V_{\rho,n}$.

    \item For every $\rho \in \Sigma(1)$ and proper subspace $W \subset V_{\rho,0}$ containing $\h$, the quotient $T/W$ of $T/H$ embeds into a toric supervariety with decorations $(N,\rho,\t,W,\{V_{\rho,i}+W\})$.

    \item For every point $x \in |X|$ of codimension at most 1, the maximal ideal $\m_x \subseteq \O_{X,x}$ is generated by a regular sequence $r_1, ..., r_d$ such that whenever $\theta \in (\Stab_\t x)_\1$ and $\ol{\theta \cdot r_i} = 0$, we have $\theta \cdot r_i \in (r_i)$.

    \item For every point $x \in |X|$ of codimension at most 1 and for every $\theta \in (\Stab_\t x)_\1$, $DS_\theta (\O_{X,x}/((\O_{X,x}^{\t_\0})_\1))$ is FR.
\end{enumerate}
\end{propdef}

\begin{proof}
    We defer this somewhat lengthy proof to Appendix \ref{section:ProofOfHR1}. In essence, we show that these conditions always hold when $[V_{\rho,0}, V_{\rho,0}]\neq0$, while for $[V_{\rho,0}, V_{\rho,0}]=0$ we mainly employ the fact that $\theta_j^\r = \d{\xi_j}$ for $\theta_j \in V_{\rho,0}$.
\end{proof}

\begin{remark}
The various equivalent characterization of $(HR_1)$ can be summarized as follows:

Characterizations (a) and (b) mean it is relatively straightforward to construct $X$ from its decorated fan. Meanwhile, (c) says that locally, away from a closed subset of codimension $\geq 2$, $X$ is a categorical limit of large-orbit toric supervarieties. Large-orbit toric supervarieties are those whose orbits have odd codimension no greater than their even codimension, and will be studied more in Section \ref{subsec:LargeOrbit}.

Characterization (d) is part of the reason we baptize this property $(HR_1)$, in analogy with the $(R_1)$ condition of Definition \ref{def:SerresCriteria}. In particular, we can think of $(HR_1)$ as an enhancement of $(R_1)$ in which we further demand that the odd vector fields behave nicely away from a closed subset of codimension $\geq 2$. (The ``homologically" part should be interpreted as arising from characterization (e), however.)

Finally, (e) links this property with the Duflo-Serganova functor $DS_\theta$ of \cite{DS} and \cite{DS20}, whose definition we recall here: If $\theta$ is an odd vector field acting semisimply on an affine supervariety $X = \Spec A$, we write $DS_\theta A$ for the cohomology of $\theta$ on the $\theta^2$-invariants of $A$, i.e.\
\begin{align*}
DS_\theta A &= \ker(\theta : A^{\theta^2}) / \Im(\theta : A^{\theta^2}) \\
&= A^\theta / \theta(A^{\theta^2}),
\end{align*}
which is naturally a superalgebra.

The quotient $\O_{X,x}/((\O_{X,x}^{\t_\0})_\1) = \O_{X,x} / ((\O_{X,x}^{\Stab_\t x})_\1)$ is the local ring at the point $x$ of a particular closed subvariety of $\Spec A_\rho$ which admits an action of the abelian supertorus $(T_\0, \t_\0 \oplus V_{\rho,0})$. Condition (e) then says that the cohomology of any odd vector field induced by this action is fermionically regular.
\end{remark}

In Example \ref{ex:WildTSV}, it is routine to check that $\Spec A$ satisfies $(HR_1)$, but $\Spec A'$ does not. In particular, consider the vector field $\theta_1$, and write $\p,\p'$ for the height-1 prime ideals of $A,A'$, respectively, generated by $t$ and the nilpotents. Then
\begin{align*}
    DS_{\theta_1} A_{\p} &= \CC[t, t\xi_2, t\xi_3, t^2 \xi_4]_{\p \cap \CC[t, t\xi_2, t\xi_3, t^2 \xi_4]} / (t) \\
    &\cong \CC[t \xi_2, t \xi_3, t^2\xi_4]
\end{align*}
is a FR superalgebra, but
\begin{align*}
    DS_{\theta_1} A'_{\p'} &= \CC[t, t\xi_2, t\xi_3, t^3 \xi_4]_{\p' \cap \CC[t, t\xi_2, t\xi_3, t^3 \xi_4]} / (t, t^2\xi_2\xi_3, t^3 \xi_4) \\
    &\cong \CC[t \xi_2, t\xi_3] / (t^2 \xi_2 \xi_3)
\end{align*}
is not.

\subsection{Admissible decorated fans}

Thanks to this means of constructing $(HR_1)$ toric supervarieties, we can now determine the exact conditions under which a decorated fan describes a $(HR_1)$ toric supervariety. We begin with the following lemma, wherein we write $V_{\sigma,0}$ for $\h_{\sigma,0}$, or equivalently $\sum_{\rho \in \sigma(1)} V_{\rho,0}$.

\begin{lemma}\label{lemma:LETTSVcontainsInducedIffDJProp}
    Let $X = \Spec A_\sigma$ be a $(HR_1)$ toric supervariety, and $W \subseteq V_{\sigma,0}$. Then $A_\sigma \supseteq L_{\CC[T]^W}(m) \neq 0$ if and only if $W$ satisfies
    \begin{enumerate}
        \item $\langle m, [W, W] \rangle  = 0$
        \item $V_{\rho,m} \subseteq W$ for all $\rho \in \sigma(1)$
        \item $\sum_{i=0}^{\langle m, u_\rho \rangle - 1} (\dim V_{\rho,i}  - \dim V_{\rho,i} \cap W) \leq \langle m, u_\rho \rangle$ for all $\rho \in \sigma(1)$
    \end{enumerate}
\end{lemma}

\begin{proof}
    First observe that condition (a) is equivalent to $L_{\CC[T]^W}(m) \neq 0$. It remains to show equivalence of (b) and (c) with $A_\sigma \supseteq L_{\CC[T]^W}(m)$.

    We begin with the case of $\sigma=\rho$ a ray. Consider the expression
    $$A_\rho = \CC[t_1, t_2^{\pm 1}, ..., t_p^{\pm 1}, t_1^{\ell_1}\xi_1, ..., t_1^{\ell_q} \xi_q]$$
    obtained in Proposition/Definition \ref{propdef:HR1}. If $W \subseteq V_{\rho,0}$ satisfies the three properties, assume without loss of generality that $\xi_1, ..., \xi_d$ is a basis for $(\t_\1/W)^*$ compatible with $$V_{\rho,0} = V_{\rho,0} + W \supseteq V_{\rho,1} + W \supseteq ... \supseteq V_{\rho,m} +W = W$$ in the sense that each $(\t_\1/(V_{\rho,j}+W))^*$ is spanned by a subset of this basis. Then $\ell_j$, for $j=1, ..., d$, is the smallest nonnegative integer such that $\xi_j$ does not vanish on $V_{\rho,j} + W$. Hence
    \begin{align*}
        \sum_{j=1}^d \ell_i &= \sum_{\ell=1}^{\langle m, u_\rho \rangle} \ell \cdot (\dim(V_{\rho, \ell-1}+W) - \dim (V_{\rho,\ell}+W)) \\
        &= \sum_{i=0}^{\langle m, u_\rho \rangle - 1} (\dim (V_{\rho,i} + W)  - \dim (V_{\rho,m} + W)) \\
        &= \sum_{i=0}^{\langle m, u_\rho \rangle - 1} (\dim V_{\rho,i} + \dim W - \dim(V_{\rho,i} \cap W) - \dim W) \\
        &= \sum_{i=0}^{\langle m, u_\rho \rangle - 1} (\dim V_{\rho,i} - \dim V_{\rho,i} \cap W) \\
        &\leq \langle m, u_\rho \rangle
    \end{align*}
    by property (c), so indeed $A_\rho$ contains all $s^m \xi_I$ for $I \subseteq \{1, ..., d\}$. Therefore $A_\sigma \supseteq L_{\CC[T]^W}(m)$.

    Conversely, if $A_\rho \supseteq L_{\CC[T]^W}(m)$, then the same chain of equalities above shows that $\sum_{i=0}^{\langle m, u_\rho \rangle - 1} (\dim V_{\rho,i} - \dim V_{\rho,i} \cap W) \leq \langle m, u_\rho \rangle$ and $V_{\rho,m} \subseteq W$.

    We now observe that properties (b) and (c) for $W \subseteq V_{\sigma,0}$ are equivalent to the same two properties for $W \cap V_{\rho,0} \subseteq V_{\rho,0}$ for each $\rho \in \sigma(1)$. Hence, if $W \subseteq V_{\sigma,0}$ satisfies (b) and (c), then
    \begin{align*}
        L_{\CC[T]^{W}}(m) &\subseteq L_{\CC[T]^{W \cap V_{\rho,0}}}(m) \\
        &\subseteq A_{\rho}
    \end{align*}
    for all $\rho$, so indeed $L_{\CC[T]^{W}}(m) \subseteq A_\sigma$. Conversely, suppose $L_{\CC[T]^{W}}(m) \subseteq A_\sigma$. Since $L_{\CC[T]^{V_{\rho,0}}}(0) \subseteq A_\rho$, it follows that $L_{\CC[T]^{W \cap V_{\rho,0}}}(m) \subseteq A_\rho$, so $W \cap V_{\rho,0}$ satisfies (b) and (c). Therefore the lemma is proven.
\end{proof}

In view of both the prior lemma and those forthcoming, we make the following definition.
\begin{definition}
    We say $W \subseteq V_{\sigma, 0}$ satisfies the \textit{dimension jumping property with respect to $\sigma$ and $m$}, or $DJ(\sigma,m)$, if $W$ satisfies the three properties of Lemma \ref{lemma:LETTSVcontainsInducedIffDJProp}. If $\rho$ is a ray, sometimes we will write $DJ(\rho,n)$ when $\langle m, u_\rho \rangle = n$.
\end{definition}

The following lemma is essentially a stronger version of Lemma \ref{lemma:DecorationCompatibility}(c). Since $W \subseteq V_{\sigma,0}$ means $[W,W] \subseteq \t_{N_\sigma}$, we note that condition (a) of Lemma \ref{lemma:LETTSVcontainsInducedIffDJProp} is equivalent to $[W,W] \subseteq \t_{N_{\sigma_m}}$.

\begin{lemma}\label{lemma:LETextistenceOfDJSubspace}
    Let $X$ be a $(HR_1)$ toric supervariety with decorated fan $(N,\Sigma,\t,\h, \{V_{\rho,i}\})$. Then for each $\sigma \in \Sigma$ and $m \in \sigma^\vee \cap M$, there is a subspace $W \subseteq V_{\sigma,0}$ satisfying $DJ(\sigma,m)$.
\end{lemma}

\begin{proof}
    Let $\sigma \in \Sigma$ and $m \in \sigma^\vee \cap M$. Then $L_{A_\sigma}(m) \neq 0$, so it contains a simple $T$-submodule $C(m)$. We have $C(m) = L_{\CC[T/S]}(m)$ for some $S \subseteq T$ with $S_\0 = T_{\ker m}$. Let $W = \s_\1 \cap V_{\sigma,0}$, so that $L_{\CC[T]^W} \subseteq L_{A_\sigma}(m)$ and hence $W$ satisfies $DJ(\sigma,m)$ by Lemma \ref{lemma:LETTSVcontainsInducedIffDJProp}.
\end{proof}

It turns out that the condition of Lemma \ref{lemma:DecorationCompatibility}(d) is sufficient to provide a converse to Lemma \ref{lemma:LETextistenceOfDJSubspace}. Therefore we will hereafter refer exclusively to the following enhanced notion of a decorated fan.

\begin{definition}\label{def:LETDecoratedFan}
    A \textit{decorated fan} is a tuple $(N,\Sigma,\t,\h,\{V_{\rho,i}\}_{\rho \in \Sigma(1), i \geq 0})$ consisting of a finite-dimensional lattice $N$, a fan $\Sigma$ in $N_\RR$, a quasiabelian Lie superalgebra $\t$ with $\t_\0 = \t_N$, an odd abelian subalgebra $\h \subseteq \t$, and a collection of descending chains $V_{\rho,0} \supseteq V_{\rho,1} \supseteq ...$ of subspaces of $\t_\1$, subject to the following compatibility conditions:
    \begin{enumerate}
        \item For each $\rho \in \Sigma(1)$, the chain $V_{\rho, 0} \supseteq V_{\rho,1} \supseteq ...$ stabilizes at $\h$.
        \item For each $\sigma \in \Sigma$ and $m \in \sigma^\vee \cap M$, there is a subspace $W \subseteq V_{\sigma,0}$ satisfying $DJ(\sigma,m)$.
        \item For each $\rho \in \Sigma(1)$, either $[V_{\rho,0}, V_{\rho,0}] = 0$ or $\codim(\h,V_{\rho,0}) = 1$.
    \end{enumerate}
\end{definition}

\begin{usethmcounterof}{thm:HR1TSVandDecoratedFanBijection}
    Any decorated fan induces a unique $(HR_1)$ toric supervariety, up to isomorphism.
\end{usethmcounterof}
\begin{proof}
    Proposition/Definition \ref{propdef:HR1} gives us both uniqueness and a method to construct the candidate toric supervariety. This method depends on conditions (b), (c), and (d) of Lemma \ref{lemma:DecorationCompatibility}, which are guaranteed by parts (a), (b), and (c), respectively, of Definition \ref{def:LETDecoratedFan}. So, for $\rho \in \Sigma(1)$, let $A_\rho$ be as in Proposition/Definition \ref{propdef:HR1}, and let $A_\sigma = \bigcap_{\rho \in \sigma(1)} A_\rho$.

    We first verify that $A_\sigma$ is finitely generated. First note that since each weight space is finite-dimensional, we may assume it contains only finitely many generators. Now choose $m \in \sigma^\vee \cap M$ such that $L_{A_\sigma}(m) = L_{\CC[T]^H}(m)$. Then we may assume every generator occurs in $L_{A_\sigma}(m')$ for $m' \leq_\sigma m$, of which there are finitely many (up to $\sigma^\perp$). We conclude that $A_\sigma$ is finitely generated, and hence that $\Spec A_\sigma$ is a (not necessarily normal) toric supervariety by Lemma \ref{ATVCharacterizationLemma}.
    
    For normality, we show that $\ol{A_\sigma} = \CC[S_\sigma]$. Since every nonzero submodule of $\CC[T]$ contains a non-nilpotent element, it suffices to show that $L_{A_\sigma}(m)$ is nonzero for every $m \in \sigma^\vee \cap M$. Using condition (b) of the definition, we obtain from Lemma \ref{lemma:LETTSVcontainsInducedIffDJProp} the desired result.

    We now claim that $A_\sigma$ localizes to $A_\rho$. It suffices to show that for each $n \geq 0$, there is $m \in \sigma^\vee \cap M$ such that $\langle m, u_\rho \rangle = n$ and $L_{A_\rho}(m) = L_{A_\sigma}(m)$. This is merely a consequence of the convex geometry of cones. In particular, the intersection $(n+\rho^\perp) \cap \sigma^\vee \cap M$ is infinite, so its elements have arbitrarily large pairings with $u_{\rho'}$ for $\rho' \in \sigma(1) \backslash \{\rho\}$. Therefore there is $m \in (n+\rho^\perp) \cap \sigma^\vee \cap M$ such that $L_{A_{\rho'}}(m) = L_{\CC[T]^H}(m)$ for all such $\rho'$, so indeed $L_{A_\rho}(m) = L_{A_\sigma}(m)$. It then follows by Proposition \ref{NormalATSVProp} that $\Spec A_\sigma$ is normal, and by Proposition/Definition \ref{propdef:HR1} that $\Spec A_\sigma$ is $(HR_1)$, completing the proof.
\end{proof}

\subsection{The abelian case}

Suppose $T$ is an abelian supertorus, and let $X$ be a toric supervariety for $T$. Since any subgroup of $T$ is normal, we may assume that the morphism $T \to X$ is an open immersion, i.e.\ $H=1$.

In this situation, since the bracket is 0, every $V_{\sigma,0}$ satisfies $DJ(\sigma,m)$, and conditions (b) and (c) of Definition \ref{def:LETDecoratedFan} are immaterial. It follows that a decorated fan can be expressed simply as a fan $\Sigma$ such that each ray $\rho \in \Sigma(1)$ is decorated by a decreasing chain of subspaces $V_{\rho,0} \supseteq V_{\rho,1} \supseteq V_{\rho,2} \supseteq ...$ of $\t_\1$ that stabilizes at 0.

\begin{example}
    We may write projective superspace $\PP^{p|q} = \Proj \CC[z_0, ..., z_p, \zeta_1, ..., \zeta_q]$ as a toric supervariety for an abelian supertorus $T$ of dimension $p|q$ as follows: As in the usual case, let $T_\0$ act by $m \in \ZZ^p$ on the monomial $z_0^r z_1^{m_1} \cdots z_p^{m_p}$ for any $r \in \ZZ$. Likewise, let $\t_\1$ act by translations on the odd part so that $\theta_i$ acts by $\d{\zeta_i}$ on the coordinate superalgebras.

    As usual, we obtain a complete fan $\Sigma$ whose rays are $\rho_i = \RR_+ e_i$ for $i=1, ..., n$ and $\rho_0 = \RR_+(-e_1-...-e_p)$. If $\sigma$ is the positive orthant $\RR_+(e_1, ..., e_p)$, then $A_\sigma = \CC[\frac{z_1}{z_0}, ..., \frac{z_p}{z_0}, \frac{\zeta_1}{z_0}, ..., \frac{\zeta_q}{z_0}] \cong \CC[t_1, ..., t_p, \xi_1, ..., \xi_q]$ and the other coordinate superalgebras of maximal cones are $\CC[t_i^{-1}t_1, ..., t_i^{-1} t_p, t_i^{-1} \xi_1, ..., t_i^{-1} \xi_q]$ for $i=1, ..., p$. Hence the decorations consist of $0 \supseteq 0 \supseteq ...$ on $\rho_1, ..., \rho_p$ and $\t_\1 \supseteq 0 \supseteq ...$ on $\rho_0$.

    Alternatively, by adjusting the characters by which $T_\0$ acts on the odd parts of the affine charts, one could arrange for any collection of subspaces $V_{\rho_i,0}$ as the decorations, so long as they are linearly independent and span $\t_\1$.
\end{example}

Since the action of $\t_\1$ on $\CC[T]$ is by $\d{\xi}$ for $\xi \in \t_\1^*$, it follows that there exists a natural splitting $\ol{A} \to A$ of the projection $A \to \ol{A}$. These splittings assemble into a global splitting $X \to X_\0$, so that any $(HR_1)$ toric supervariety for an abelian supertorus is \textit{split} in this sense.

\subsection{The case of large orbits}\label{subsec:LargeOrbit}
In this section we will assume that $\codim(\h, V_{\rho,0}) \leq 1$, so that the odd codimension of an orbit never exceeds the even codimension. Such a toric supervariety will be said to have \textit{large orbits}. In view of Proposition/Definition \ref{propdef:HR1}, it is clear that toric supervarieties with large orbits are always $(HR_1)$. Moreover, $\sum_{\rho \in \sigma(1)} V_{\rho,m}$ always satisfies $DJ(\sigma,m)$, so condition (b) of Definition \ref{def:LETDecoratedFan} can be weakened to the condition of Lemma \ref{lemma:DecorationCompatibility}(c). These observations are summarized in the following definition and proposition:

\begin{definition}\label{def:LargeOrbitDecoratedFan}
    A \textit{large-orbit decorated fan} is a tuple $(N,\Sigma,\t,\h,\{V_{\rho,i}\})$ consisting of a finite-dimensional lattice $N$, a fan $\Sigma$ in $N_\RR$, a quasiabelian Lie superalgebra $\t$ with $\t_\0 = \t_N$, an odd abelian subalgebra $\h \subseteq \t$, and a collection of descending chains $V_{\rho,0} \supseteq V_{\rho,1} \supseteq ...$ of subspaces of $\t_\1$, subject to the following compatibility conditions:
    \begin{enumerate}
        \item For each $\rho \in \Sigma(1)$, the chain $V_{\rho, 0} \supseteq V_{\rho,1} \supseteq ...$ stabilizes at $\h$.
        \item If $\rho,\rho' \in \sigma(1)$ and $m \in \sigma^\vee \cap M$, then $[V_{\rho,m}, V_{\rho',m}] \subseteq \t_{N_{\sigma_m}}$.
        \item For each $\rho \in \Sigma(1)$, $\codim(\h,V_{\rho,0}) \leq 1$.
    \end{enumerate}
\end{definition}

In particular, each chain $V_{\rho,0} \supseteq V_{\rho,1} \supseteq ...$ is constant except possibly for a one-dimensional jump at a single index.

\begin{proposition}\label{prop:largeOrbitClassification}
    Up to $T$-equivariant isomorphism, there is a bijective correspondence between large-orbit decorated fans and toric supervarieties with large orbits.
\end{proposition}

For supertori such as $T=Q(1)^n$ that have no nonzero square roots of 0 (or more generally, no purely odd subalgebras properly containing $\h$), every toric supervariety has large orbits by Lemma \ref{lemma:DecorationCompatibility}(d). Moreover, since $[V_{\rho,1}, V_{\rho,1}]=0$, it holds that $V_{\rho,1}=\h$. Consequently, all the data of such a toric supervariety are concentrated in the fan decorations $V_{\rho,0}$ for $\rho \in \Sigma(1)$. The data of a decorated fan can then be simplified as in the following corollary.

\begin{corollary}\label{cor:noNonzeroSquareRootsOfZero}
    Let $T$ be a supertorus such as $Q(1)^n$ with no nonzero square roots of 0. Then a decorated fan with $\t=\Lie T$ can be reduced to the data of $(N,\Sigma,\t, \{V_\rho\}_{\rho \in \Sigma(1)})$ satisfying the following compatibility conditions.
    \begin{enumerate}
        \item If $\rho,\rho' \in \sigma(1)$, then $[V_\rho, V_{\rho'}] \subseteq \t_{N_\sigma}$.
        \item For each $\rho \in \Sigma(1)$, $\dim(V_\rho)\leq 1$.
    \end{enumerate}
\end{corollary}

We will explore some examples of such toric supervarieties presently, as well as in section \ref{QGrSection}.

\begin{example}\label{ex:PnFaithfulQ(1)n}
    Let us use Proposition \ref{prop:largeOrbitClassification} and Corollary \ref{cor:noNonzeroSquareRootsOfZero} to classify toric supervarieties with supertorus $Q(1)^n$ and underlying variety $\PP^n \cong X_\Sigma$ for the complete fan whose rays are $\rho_i = \RR_+x_i$ for $i=1, ..., n$ and $\rho_0 = \RR_+(-x_1-...-x_n)$.

    \begin{figure}
        \centering
        \begin{tikzpicture}[scale=2]
        
            \draw[thick, ->] (0,0) -- (1,0) node[right] {};
            \draw[thick, ->] (0,0) -- (0,1) node[above] {};
            \draw[thick, ->] (0,0) -- (-0.7,-0.7) node[below left] {};
            
            \node at (1.6, 0) {$V_{\rho_1} = \CC \theta_1$};
            \node at (0, 1.2) {$V_{\rho_2} = \CC \theta_2$};
            \node at (-0.84, -0.95) {$V_{\rho_0} = 0$};
        
            \fill[orange, opacity=0.3] (0,0) -- (1,0) -- (-0.7,-0.7) -- cycle;
            \fill[blue, opacity=0.3] (0,0) -- (1,0) -- (0,1) -- cycle;
            \fill[green, opacity=0.3] (0,0) -- (0,1) -- (-0.7,-0.7) -- cycle;
        
        \end{tikzpicture}
        \caption{The decorated fan of an action of $Q(1)^2$ on $\PP^{2|2}$}
        \label{fig:P2decoratedFan}
    \end{figure}
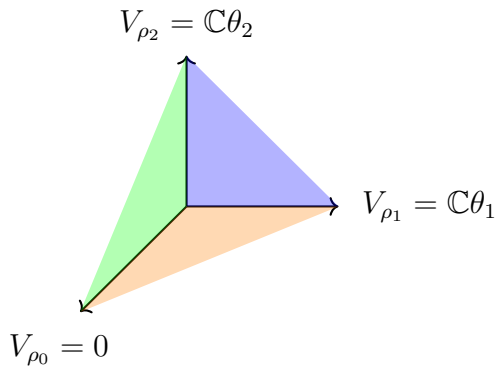

    For $i=1, ..., n$, the ray $\rho_i$ must be decorated by a subspace $V_{\rho_i}$ such that $[V_{\rho_i}, V_{\rho_i}] \subseteq \CC x_i$. Hence $V_{\rho_i} = \CC \theta_i$ or 0. Likewise, unless $V_{\rho_0}=0$, we have $V_{\rho_0} = \CC(\theta_1 \pm \theta_2 \pm ... \pm \theta_n)$ for some $2^{n-1}$ choices of $\pm$.

    It is straightforward to verify that condition (b) of the definition of a large-orbit decorated fan holds regardless of which subspaces are chosen. We therefore obtain a collection of $2^n(1+2^{n-1})$ toric supervarieties which are not equivariantly isomorphic. Many, however, are isomorphic via toric morphisms (to be defined in the following section).

    Notice that if all but one of these decorations is nonzero, then the decorated fan describes a supervariety isomorphic to projective superspace $\PP^{n|n}$. For instance, if $V_{\rho_i} = \CC \theta_i$ and $V_{\rho_0} = 0$, then the coordinate superalgebras of the affine charts can be written as $$\CC[t_1, ..., t_n, t_1\xi_1, ..., t_n \xi_n]$$ and $$\CC[t_i^{-1} t_1, ..., t_i^{-1}, ..., t_i^{-1} t_n \xi_1, ..., \xi_i, ..., t_i^{-1} t_n \xi_n].$$ Figure \ref{fig:P2decoratedFan} depicts the corresponding decorated fan for $n=2$,
    
    If instead we change the decoration of $\rho_0$ to $\theta_1 + ... + \theta_n$, then the resulting affine charts have coordinate superalgebras
    $$\CC[t_1, ..., t_n, t_1\xi_1, ..., t_n \xi_n]$$ and $$\CC[t_i^{-1} t_1(1+\xi_i\xi_1), ..., t_i^{-1}, ..., t_i^{-1} t_n(1+\xi_i \xi_n), t_i^{-1} t_1(\xi_i-\xi_1), ..., t_i^{-1} \xi_i, ..., t_i^{-1} t_n (\xi_i-\xi_n)],$$
    so the supervariety is decidedly not isomorphic to projective superspace. Figure \ref{fig:superPolytopeAndFan}(B) depicts corresponding decorated fan for $n=1$.
\end{example}

In general, when $T=Q(1)^n$, there are finitely many possible decorations for each ray. Namely, for $\rho = \RR_+ (a_1x_1 + ... + a_nx_n)$, there are $2^{d-1}$ many possible ``square root subspaces" $\CC (\sqrt{a_1} \theta_1 \pm ... \pm \sqrt{a_n} \theta_n)$, where $d$ is the number of indices $i=1, ..., n$ for which $a_i \neq 0$.

The prior example admitted no issues of compatibility between different rays of the same cone. This is not ordinarily the case; if $\rho_1 = \RR_+(a_1x_1 + ... + a_nx_n)$ and $\rho_2 = \RR_+(b_1x_1+...+b_nx_n)$, compatibility is less common if the linear matroid on the $n$ vectors
$$\begin{pmatrix}
    a_1 \\ b_1
\end{pmatrix}, ..., \begin{pmatrix}
    a_n \\ b_n
\end{pmatrix}$$
has certain independence properties. For example, if two rays span a 2-dimensional cone and the corresponding matroid contains a 3-circuit, then one of the decorations $V_\rho$ must be 0. While we do not explore this idea in depth in this paper, we provide the following example of this phenomenon.

\begin{example}
    Let $n=3$ and consider the cone whose rays are $\rho_1 = \RR_+(x_1+x_2)$ and $\rho_2 = \RR_+(x_1+x_3)$. If $V_{\rho_1}$ and $V_{\rho_2}$ are both nonzero, then $[V_{\rho_1}, V_{\rho_2}] = \CC x_1$, which is not contained within the span of $\rho_1$ and $\rho_2$. Hence the decorations are compatible if and only if at least one of them is 0.
\end{example}

\subsection{Decorated polytopes}
Classically, a toric variety $X_\Sigma$ is projective if and only if $\Sigma$ is the normal fan $\Sigma_P$ of a lattice polytope $P$ in $M_\RR := M \otimes_\ZZ \RR$. In particular, a face $F$ of dimension $r$ induces a cone $\sigma_F$ of codimension $r$, which in turn corresponds to a torus orbit of dimension $r$. 

Suppose $X$ is a toric supervariety whose underlying variety is $X_{\Sigma_P}$. Due to the bijection between faces of $P$ and cones of $\Sigma_P$, we may assign the same data $V_{\sigma_F, 0} \supseteq V_{\sigma_F,1} \supseteq ...$ to each face $F$ of codimension 1 in $P$. Since $M_\RR$ is a real form of $\t_\0^*$ and $V_{\sigma_F} \subseteq \t_\1$, we will prefer the assignment $W_{F,i} := (\t_\1/V_{\sigma_F,i})^* \subseteq \t_\1^*$ when attempting to decorate such a polytope. 

Of course, in order for a decorated polytope to define a toric supervariety, one requires compatibility conditions on the decorations $W_{F,i}$. We neglect to state these on the grounds that they are essentially the same as those in Definitions \ref{def:LETDecoratedFan} and \ref{def:LargeOrbitDecoratedFan}, and that they are easier to express in terms of the $V_{\rho,i}$.

Suppose the decorated polytope $(P, \{W_{F,i}\})$ corresponds to a toric supervariety with large orbits, and define $W_F = (\t_\1/V_{\sigma_F, 0})^*$ for all faces $F \subseteq P$. We may then view the decorated polytope as the subset $\P = \bigcup_{F \in P} \Int F \times W_F$ of $M_\RR \oplus \t_\1^*$, where $\Int F$ is the relative interior of $F$. This subset fails to be bounded or even a polyhedron, but it is at least convex and locally closed. The upshot of this viewpoint will be explored more in section \ref{QGrSection}.

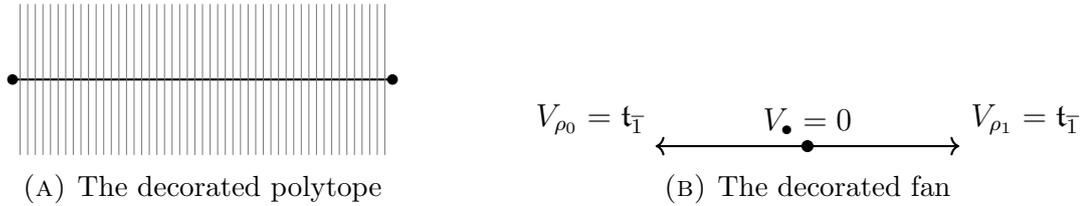
\begin{figure}
        \centering
        
        \begin{subfigure}{.5\textwidth}
            \centering
            \begin{tikzpicture}
                \draw[thick] (0,0) -- (5,0);
                
                \foreach \x in {0.1,0.2,...,4.9} {
                    \draw[black!50] (\x,0) -- ++(0,1); 
                    \draw[black!50] (\x,0) -- ++(0,-1);
                }
        
                \fill[black] (0,0) circle (0.07);
                \fill[black] (5,0) circle (0.07);
                
            \end{tikzpicture}
            \caption{The decorated polytope}
            \label{fig:enter-label}
        \end{subfigure}%
        \begin{subfigure}{.5\textwidth}
          \centering
          \begin{tikzpicture}[scale=2]
                \draw[thick, ->] (0,0) -- (1,0) node[above right] {$V_{\rho_1} = \t_\1$};
                \draw[thick, ->] (0,0) -- (-1,0) node[above left] {$V_{\rho_0} = \t_\1$};
            
                \fill[black] (0,0) circle (0.04) node [above] {$V_\bullet = 0$};
            \end{tikzpicture}
          \caption{The decorated fan}
        \end{subfigure}
        
        \caption{Two different combinatorial descriptions of $\QGr(1,2)$ as a toric supervariety for $Q(1)$.}
        \label{fig:superPolytopeAndFan}
    \end{figure}

\begin{example}
    Continuing Example \ref{ex:PnFaithfulQ(1)n} in the case of $n=1$, we note that the decorations $V_{\rho_0} = V_{\rho_1} = \t_\1$ describe the ``isomeric projective line," or $\QGr(1,2)$, to be defined in section \ref{QGrSection}. The corresponding decorated polytope and decorated fan are illustrated in Figure \ref{fig:superPolytopeAndFan}.
\end{example}

\begin{figure}
    
\end{figure}

\subsection{Toric morphisms}

Towards a categorical version of our classification of toric supervarieties, we now consider toric morphisms. Let $X$ and $X'$ be toric supervarieties with open orbits $T/H$ and $T'/H'$, respectively.

\begin{definition}
    A \textit{toric morphism} $X \to X'$ is a commutative diagram
    \begin{center}
    \begin{tikzcd}
        T \arrow[r] \arrow[d, "\phi"] & T/H \arrow[r], \arrow[d, "\tilde \phi"] & X \arrow[d, "\varphi"] \\
        T' \arrow[r] & T'/H' \arrow[r] & X'
    \end{tikzcd}
    \end{center}
    where $\phi$ is a morphism of supergroups such that $\phi(H) \subseteq H'$, $\tilde \phi$ is the induced map $T/H \to T'/H'$, and $\varphi : X \to X'$ extends $\tilde \phi$.
\end{definition}

A toric morphism of toric supervarieties naturally drops to a toric morphism of the underlying toric varieties. So, we obtain from the pair $(\phi, \varphi)$ a morphism $\phi^\circ : N \to N'$ of lattices such that for each $\sigma \in \Sigma$, there is $\sigma' \in \Sigma'$ for which $\phi^\circ_\RR(\sigma) \subseteq \sigma'$. Such a $\phi^\circ$ is called a \textit{morphism of fans}. 

In addition, we obtain a morphism $d\phi : \t \to \t'$ of Lie superalgebras such that $\phi^\circ_\CC = d \phi_\0$. It remains to describe the compatibility conditions with respect to the decorations. Suppose $\phi^\circ_\RR(\sigma) \subseteq \sigma'$, so we may reduce to the affine case $\varphi : \Spec A_\sigma \to \Spec A'_{\sigma'}$. The map $\varphi$ dualizes to $\varphi^* : A'_{\sigma'} \to A_\sigma$, a restriction of $\phi^* : \CC[T'] \to \CC[T]$. The underlying map of rings $\ol{\varphi^*} : \kk[\sigma'^\vee \cap M'] \to \kk[\sigma^\vee \cap M]$ sends $m' \in M' \cong \Hom(N', \ZZ)$ to $m' \circ \phi^\circ \in M \cong \Hom(N,\ZZ)$, so it must hold that $\phi^*(L_{A'_{\sigma'}}(m')) \subseteq L_{A_\sigma}(m' \circ \phi^\circ)$.

Thus, a toric morphism is determined by a morphism of fans $\phi^\circ : (N,\Sigma) \to (N',\Sigma')$ together with a compatible morphism of Lie superalgebras $d\phi : \t \to \t'$ such that $\phi^*(L_{A'_{\sigma'}}(m')) \subseteq L_{A_\sigma}(m' \circ \phi^\circ)$ whenever $\phi^\circ_\RR(\sigma) \subseteq \sigma'$ and $m' \in \sigma'^\vee \cap M$. One would hope that a condition involving only the decorations $V_{\rho,i}$ and $V'_{\rho',i}$ would suffice. For large-orbit toric supervarieties, such a condition is indeed viable due to the following property.

\begin{lemma}\label{lemma:largeOrbitSumOfInducedReps}
    Let $X$ be a toric supervariety with large orbits, and let $(N,\Sigma,\t,\h,\{V_{\rho,i}\})$ be its decorated fan. Then for $\sigma \in \Sigma$, it holds that $L_{A_\sigma}(m) = L_{\CC[T]^{V_{\sigma,m}}}(m)$ where $V_{\sigma,m} = \sum_{\rho \in \sigma(1)} V_{\rho,m}$.
\end{lemma}
\begin{proof}
    The claim is immediate in the case of $\dim \sigma \leq 1$ by construction of $A_\sigma$ in Proposition/Definition \ref{propdef:HR1}. Otherwise we have
    \begin{align*}
        L_{A_\sigma}(m) &= \bigcap_{\rho \in \sigma(1)} L_{A_\rho}(m) \\
        &= \bigcap_{\rho \in \sigma(1)} L_{\CC[T]^{V_{\rho,m}}}(m) \\
        &= L_{\CC[T]^{V_{\sigma,m}}}(m)
    \end{align*}
    and we are finished.
\end{proof}

Hence the condition that $\phi^*(L_{A'_{\sigma'}}(m')) \subseteq L_{A_\sigma}(m' \circ \phi^\circ)$ is equivalent to $\phi^*(L_{\CC[T']^{V'_{\sigma',m'}}}(m')) \subseteq L_{\CC[T]^{V_{\sigma,m' \circ \phi^\circ}}}(m' \circ \phi^\circ)$ or simply $d\phi(V_{\sigma,m' \circ \phi^\circ}) \subseteq V'_{\sigma',m'}$.

\begin{definition}\label{def:LargeOrbitDecFanMorphism}
    Let $(N,\Sigma,\t,\h,\{V_{\rho,i}\})$ and $(N',\Sigma',\t',\h',\{V'_{\rho',i}\})$ be large-orbit decorated fans, with $\phi^\circ : (N,\Sigma) \to (N',\Sigma')$ a morphism of fans and $d\phi : \t \to \t'$ a compatible morphism of Lie superalgebras. We say $d\phi$ is \textit{compatible with the decorations} if $d\phi(V_{\sigma,m' \circ \phi^\circ}) \subseteq V_{\sigma', m'}$ for any $\sigma \in \Sigma$ and $\sigma' \in \Sigma$ such that $\phi^\circ_\RR(\sigma) \subseteq \sigma'$, and any $m' \in \sigma'^\vee \cap M'$. In this situation, we call the pair $(\phi^\circ, d\phi)$ a \textit{morphism of large-orbit decorated fans}.
\end{definition}

That is, we have the following proposition:

\begin{proposition}\label{prop:LargeOrbitMorphismClassifcation}
    Let $X$ and $X'$ be toric supervarieties with large orbits. A morphism $\phi : T \to T'$ of supergroups drops to $T/H \to T'/H'$ and subsequently extends to a toric morphism $X \to X'$ if and only if $(\phi^\circ, d\phi)$ is a morphism of large-orbit decorated fans.
\end{proposition}
\begin{proof}
    As above, this follows from Lemma \ref{lemma:largeOrbitSumOfInducedReps}.
\end{proof}

We may then define the category of large-orbit decorated fans via the objects and morphisms of Definitions \ref{def:LargeOrbitDecoratedFan} and \ref{def:LargeOrbitDecFanMorphism}. Combining Propositions \ref{prop:largeOrbitClassification} and \ref{prop:LargeOrbitMorphismClassifcation}, we immediately obtain the following theorem.

\begin{usethmcounterof}{thm:largeOrbitEquivalenceOfCats}
    The category of large-orbit toric supervarieties is equivalent to the category of large-orbit decorated fans.
\end{usethmcounterof}

Generalizing this equivalence to $(HR_1)$ toric supervarieties turns out to be surprisingly difficult. While it is not difficult to show, in analogy with Lemma \ref{lemma:largeOrbitSumOfInducedReps}, that $A_\rho$ is spanned by $L_{\CC[T]^W}(m)$ for $W$ satisfying $DJ(\rho,m)$, it is not the case that $A_\sigma$ is spanned by $L_{\CC[T]^W}(m)$ for $W$ satisfying $DJ(\sigma,m)$. The following example demonstrates the strangeness of the $(HR_1)$ situation.

\begin{example}
    In this example we construct a $(HR_1)$ toric supervariety for which some $L_{A_\sigma}(m)$ is not spanned by induced representations. Consider the $(3|4)$-dimensional abelian supertorus $T = (\CC^\times)^3 \times \AA^{0|4}$. As before, we write $\CC[T] = \CC[t_1^{\pm 1}, t_2^{\pm 1}, t_3^{\pm 1}, \xi_1, \xi_2, \xi_3, \xi_4]$ and $\{\theta_1, \theta_2, \theta_3, \theta_4\}$ a dual basis of $\t_\1$ to $\{\xi_1, \xi_2, \xi_3, \xi_4\}$, so that $\theta_i$ acts by $\d{\xi_i}$ on $\CC[T]$.
    
    Let $\sigma$ be the positive orthant in $\RR^3$, with rays $\rho_i = \RR_+ x_i$ for $i=1, 2, 3$. We decorate the rays as follows:
    \begin{center}
    \begin{tikzcd}[column sep=small]
        V_{\rho_1,0} \arrow[r, phantom, sloped, "\supseteq"]\arrow[d, phantom, sloped, "="] & V_{\rho_1,1} \arrow[r, phantom, sloped, "\supseteq"]\arrow[d, phantom, sloped, "="] & V_{\rho_1,2} \arrow[r, phantom, sloped, "\supseteq"]\arrow[d, phantom, sloped, "="] & V_{\rho_1,3} \arrow[r, phantom, sloped, "\supseteq"]\arrow[d, phantom, sloped, "="] & V_{\rho_1,4} \arrow[d, phantom, sloped, "="] \\
        \CC\{\theta_1, \theta_2, \theta_3, \theta_4\} & \CC\{\theta_2, \theta_3, \theta_4\} & \CC\{\theta_2, \theta_3\} & \CC\{\theta_2\} & 0
    \end{tikzcd}
    \end{center}
    and
    \begin{center}
    \begin{tikzcd}[column sep=small]
        V_{\rho_2,0} \arrow[r, phantom, sloped, "\supseteq"]\arrow[d, phantom, sloped, "="] & V_{\rho_2,1} \arrow[r, phantom, sloped, "\supseteq"]\arrow[d, phantom, sloped, "="] & V_{\rho_2,2} \arrow[r, phantom, sloped, "\supseteq"]\arrow[d, phantom, sloped, "="] & V_{\rho_2,3} \arrow[r, phantom, sloped, "\supseteq"]\arrow[d, phantom, sloped, "="] & V_{\rho_2,4} \arrow[d, phantom, sloped, "="] \\
        \CC\{\theta_1, \theta_2, \theta_3, \theta_4\} & \CC\{2\theta_1 + \theta_3, \theta_2, \theta_4\} & \CC\{2\theta_1 + \theta_3, \theta_2 - 2\theta_4\} & \CC\{\theta_2 - 2\theta_4\} & 0
    \end{tikzcd}
    \end{center}
    and
    \begin{center}
    \begin{tikzcd}[column sep=small]
        V_{\rho_3,0} \arrow[r, phantom, sloped, "\supseteq"]\arrow[d, phantom, sloped, "="] & V_{\rho_3,1} \arrow[r, phantom, sloped, "\supseteq"]\arrow[d, phantom, sloped, "="] & V_{\rho_3,2} \arrow[r, phantom, sloped, "\supseteq"]\arrow[d, phantom, sloped, "="] & V_{\rho_3,3} \arrow[r, phantom, sloped, "\supseteq"]\arrow[d, phantom, sloped, "="] & V_{\rho_3,4} \arrow[d, phantom, sloped, "="] \\
        \CC\{\theta_1, \theta_2, \theta_3, \theta_4\} & \CC\{\theta_1, \theta_2 - \theta_4, \theta_3\} & \CC\{\theta_1 + \theta_3, \theta_2 - \theta_4\} & \CC\{\theta_1 + \theta_3\} & 0
    \end{tikzcd}
    \end{center}

    Hence
    \begin{align*}
        A_{\rho_1} &= \CC[t_1, t_2^{\pm 1}, t_3^{\pm 1}, t_1 \xi_1, t_1^2 \xi_4, t_1^3 \xi_3, t_1^4 \xi_2] \\
        A_{\rho_2} &= \CC[t_1^{\pm 1}, t_2, t_3^{\pm 1}, t_2 (\xi_1-2\xi_3), t_2^2 (2\xi_2+\xi_4), t_2^3 (2\xi_1-\xi_3), t_2^4 (\xi_2+2\xi_4)] \\
        A_{\rho_3} &= \CC[t_1^{\pm 1}, t_2^{\pm 1}, t_3, t_3(\xi_2+\xi_4), t_3^2 (\xi_1-\xi_3), t_3^3 \xi_4, t_3^4 \xi_1]
    \end{align*}
    and a straightforward calculation shows that 
    \begin{align*}
    L_{A_\sigma}(5, 5, 5) &= \bigcap_{i=1}^3 L_{A_{\rho_i}}(5, 5, 5) \\
    &=\CC t_1^5t_2^5t_3^5\{1, \xi_1, \xi_2, \xi_3, \xi_4, \xi_1\xi_2 + \xi_3\xi_4\}  .  
    \end{align*}
    Since $\xi_1\xi_2+\xi_3\xi_4$ is an indecomposable element of $\bigwedge \t_\1^*$, it follows that $L_{A_\sigma}(5, 5, 5)$ is not the span of induced representations of the form $L_{\CC[T]^W}(5,5,5)$.
\end{example}

\subsection{Smoothness of toric supervarieties}

In this section, we will provide a criterion for smoothness of a toric supervariety based on its decorated fan.

In the classical situation, a toric variety $X$ is smooth if and only if each cone $\sigma$ in its fan $\Sigma$ is a \textit{smooth cone} in the sense that the primitive generators of its rays form part of a basis for $N$. This has the effect that every smooth affine toric variety is isomorphic to $\CC^n \times (\CC^\times)^{p-n}$ for some $p \geq n \geq 0$, a fact which we will use in our proof of Theorem \ref{SmoothTSVTheorem}.

\begin{definition}
    Let $X$ be a supervariety and $x \in X(\kk)$ a closed point. We say $X$ is \textit{smooth} at $x$ if $\O_{X,x}$ is a regular local superring. We say $X$ is \textit{smooth} if it is smooth at all closed points.
\end{definition}

There exist many equivalent formulations for smoothness of a supervariety (or more broadly, smoothness of a morphism of superschemes), several of which can be found in Appendix B of \cite{Sherman}.

\begin{theorem}\label{SmoothTSVTheorem}
    Let $X$ be a $(HR_1)$ toric supervariety with decorated fan $(N,\Sigma,\t,\h,\{V_{\rho,i}\})$. Then $X$ is smooth if and only if for every cone $\sigma \in \Sigma$, the following three conditions hold:
    \begin{enumerate}
        \item $\sigma$ is a smooth cone.
        \item There is a basis of $\t_\1$ such that for all $\rho \in \sigma(1)$ and $i \geq 0$, $V_{\rho,i}$ is the span of a subset of that basis.
        \item For each $\rho \in \sigma(1)$, if $\codim(\h, V_{\rho,0}) \geq 2$, then $[V_{\rho,0}, V_{\rho',0}] \subseteq \t_{N_{\rho'}}$ for all $\rho' \in \sigma(1)$.
    \end{enumerate}
\end{theorem}

\begin{proof}
    We may assume $X$ is affine, and that $X_\0 \cong \CC^n \times (\CC^\times)^{p-n}$, so $X$ corresponds to the cone $\sigma$ whose rays are $\rho_i = \RR_+ e_i$ for $i=1, ..., n$, where $e_i$ is the $i$th unit vector in $\ZZ^{p}$. If $X$ is smooth, then (a) is immediate and we may write
    $$A_\sigma = \CC[t_1, ..., t_n, t_{n+1}^{\pm 1}, ..., t_{p}^{\pm 1}, t^{m_1}(\xi_1+...), ..., t^{m_q} (\xi_q+...)]$$
    for some non-nilpotent $t_i$, linearly independent $\xi_i \in \t_\1^*$, and products of nilpotent terms represented by ellipses. Localizing to some $A_{\rho_i}$, we see that we may assume each $\xi_j + ...$ simply equals $\xi_j$. We therefore have
    $$A_\sigma = \CC[t_1, ..., t_n, t_{n+1}^{\pm 1}, ..., t_{p}^{\pm 1}, t^{m_1}\xi_1, ..., t^{m_q} \xi_q]$$
    and in particular
    $$A_{\rho_i} = \CC[t_1^{\pm 1}, ..., t_i, ..., t_p^{\pm 1}, t_i^{\langle m_1, e_i \rangle} \xi_1, ..., t_i^{\langle m_q, e_i \rangle} \xi_q]$$
    for $i=1, ..., n$. Choosing a basis of $\t_\1^*$ that includes the $\xi_i$ yields a dual basis $\theta_1, ..., \theta_q, \theta_{q+1}, ..., \theta_r$ of $\t_\1$ for which $V_{\rho_i,m} = \CC\{\theta_j \mid \langle m, e_i \rangle < \langle m_j, e_i \rangle \} \oplus \CC\{\theta_{q+1}, ..., \theta_r\}$. Hence (b) holds.

    Now for (c), note that if $\codim(\h, V_{\rho_i,0}) \geq 2$, then by construction of $A_\rho$ in Proposition/Definition \ref{propdef:HR1}, we have $t_i \in \CC[T]^{V_{\rho_i,0}}$. Since $t_i \in \CC[T]^{V_{\rho_\ell,0}}$ for all $\ell \neq i$ as well, it follows that $t_i \in \CC[T]^{[V_{\rho_i,0}, V_{\rho_\ell,0}]}$ and therefore $\langle \epsilon_i, [V_{\rho_i,0}, V_{\rho_\ell,0}] \rangle = 0$, where $\{\epsilon_i\}$ in $M$ is dual to $\{e_i\}$ in $N$. Hence $[V_{\rho_i,0}, V_{\rho_\ell,0}] \subseteq \t_{N_{\rho_\ell}}$ and we are finished.
    
    Conversely, suppose properties (a), (b), and (c) hold. As above, property (c) allows us to find $t_i \in \CC[T]^{\sum_{\ell=1}^n V_{\rho_\ell,0}}$ if $\codim(\h, V_{\rho_i,0}) \geq 2$, or otherwise in $\CC[T]^{V_{\rho_i,1} + \sum_{\ell \neq i} V_{\rho_{\ell},0}}$. Now pick such a basis $\theta_1, ..., \theta_r$ as in (b), where $\h = \CC\{\theta_{q+1}, ..., \theta_r\}$, and let $\xi_1, ..., \xi_r$ be the dual basis of $\t_\1^*$. Then each
    $$A_{\rho_i} = \CC[t_1^{\pm 1}, ..., t_i, ..., t_p^{\pm 1}, t_i^{n_{i1}} \xi_1, ..., t_i^{n_{iq}} \xi_q]$$
    where $n_{ij} \geq 0$ is minimal such that $\theta_j \notin V_{\rho_i, n_{ij}}$. Then let $m_j = n_{1j}\epsilon_1 + ... + n_{qj} \epsilon_q$, so
    $$A_\sigma = \CC[t_1, ..., t_n, t_{n+1}^{\pm 1}, ..., t_{p}^{\pm 1}, t^{m_1}\xi_1, ..., t^{m_q} \xi_q]$$
    and indeed $X$ is smooth.
\end{proof}

This theorem simplifies in the case of a toric supervariety with large orbits. In particular, condition (c) is immaterial, so smoothness in the odd directions depends entirely on finding such a basis as in (b).

This condition on basis elements is equivalent to Klyachko's condition \cite{Klyachko} that a descending filtration $E = ... \supseteq E^\rho(-1) \supseteq E^\rho(0) \supseteq E^\rho(1) \supseteq ... = 0$ for each $\rho \in \Sigma(1)$ specifies a toric vector bundle with fiber $E$ over the identity. To obtain a Klyachko filtration from our $V_{\rho,0} \supseteq V_{\rho,1} \supseteq ...$, one must take $E^\rho(i) = 0$ for $i>0$ and $E^\rho(i) = (\t_\1/V_{\rho,-i})^*$ for $i \leq 0$. These Klyachko data yield the vector bundle $\J_X / \J_X^2$, the conormal bundle of $X_\0$ in $X$.

Note the sign discrepancy in $E^\rho(i) = (\t_\1/V_{\rho,-i})^*$; this occurs because $E^\rho(i)$ is the $i m_1$-eigenspace of $J_{A_\rho} / J_{A_\rho}^2$ for the left regular representation, wherein the even torus acts by the character $-m$ on the monomial $t^m$.

Using this characterization of smoothness, it is not difficult to construct a toric resolution of singularities for a toric supervariety.

\begin{proposition}\label{prop:existenceOfResolutionOfSingularities}
    Let $X$ be a toric supervariety for the supertorus $T$. Then there exists a smooth toric supervariety $X'$ for $T$ and a proper, birational, $T$-equivariant toric morphism $X' \to X$.
\end{proposition}

\begin{proof}
    Let $(N,\Sigma,\t,\h,\{V_{\rho,i}\})$ be the decorated fan of $X$. Construct a refinement $\Sigma'$ of $\Sigma$ (with the same support) by introducing new rays until each maximal cone is smooth and contains at most one ray from the original fan $\Sigma$. Decorate each new ray $\rho'$ by $V'_{\rho',i} = \h$ for all $i \geq 0$. Since each cone contains at most one of the original rays, the resulting decorated fan $(N, \Sigma', \t, \h, \{V'_{\rho',i}\})$ satisfies the conditions of Theorem \ref{SmoothTSVTheorem}.

    We now verify that the identity $T \to T$ extends to $X' \to X$. It suffices to consider a ray $\rho' \in \Sigma'$ and a cone $\sigma \in \Sigma$ such that $\rho' \subseteq \sigma$. If $\rho' \notin \Sigma$, then $\rho'$ has trivial decorations and so each $L_{A_{\rho'}}(m) = L_{\CC[T]^H}(m) \supseteq L_{A_\sigma}(m)$. Otherwise $\rho' \in \Sigma$, in which case $A_{\rho'} = A_\rho$. Therefore we obtain a $T$-equivariant morphism $X' \to X$. Birationality is immediate since the open orbits are both isomorphic to $T/H$. Finally, properness of the underlying morphism of varieties follows since the fans have the same support; this is sufficient since properness depends only on the underlying varieties (see e.g.\ \cite{BRP}).
\end{proof}

\begin{example}\label{ex:SingularButEvenSmooth}
    The toric supervariety $$X = \Spec \CC[t_1, t_2, t_3, t_1t_2 \xi_1, t_1t_3 \xi_2, t_2t_3(\xi_1+\xi_2), t_1t_2t_3 \xi_1\xi_2]$$
    has smooth underlying variety, but is not itself smooth. Its failure to be smooth can be detected by checking that the subspaces $\CC \theta_2, \CC \theta_1$, and $\CC(\theta_1-\theta_2)$ occur among the decorations. The decorated fan for $X$, as well as that of a (non-affine) resolution of singularities $X'$, can be found in Figure \ref{fig:ResolutionOfSingularities}. Note that the resolution depicted is not the one constructed in Proposition \ref{prop:existenceOfResolutionOfSingularities}, since the fan is not sufficiently refined.

    \tdplotsetmaincoords{60}{110}
    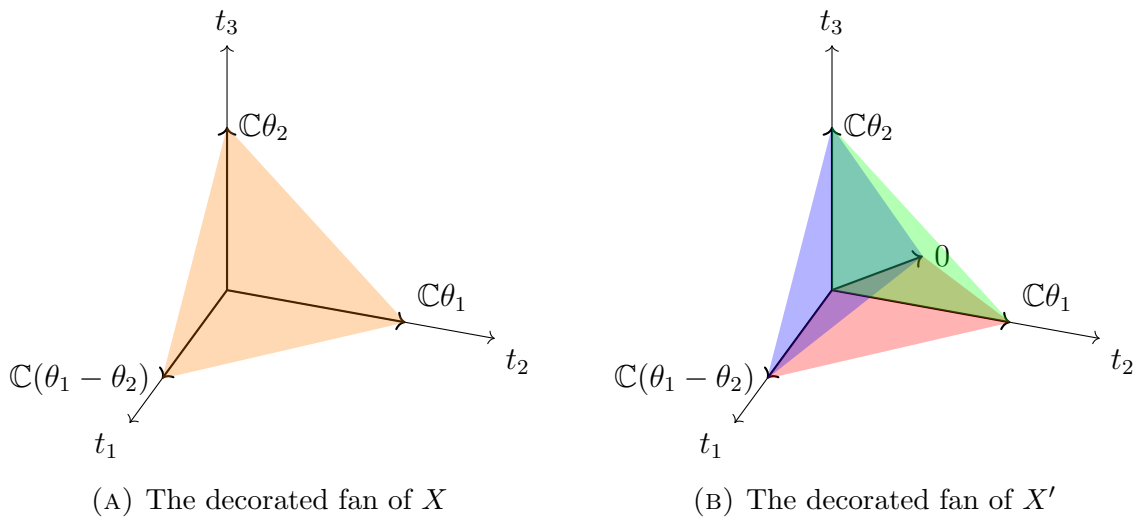
\begin{figure}
        \centering
        \begin{subfigure}{.5\textwidth}
          \centering
            \begin{tikzpicture}[tdplot_main_coords, scale=2.5]
                \draw[->] (0,0,0) -- (1.5,0,0) node[anchor=north east] {$t_1$};
                \draw[->] (0,0,0) -- (0,1.5,0) node[anchor=north west] {$t_2$};
                \draw[->] (0,0,0) -- (0,0,1.5) node[anchor=south] {$t_3$};
            
                \draw[thick, ->] (0,0,0) -- (1,0,0) node[anchor=east] {$\CC(\theta_1-\theta_2)$};
                \draw[thick, ->] (0,0,0) -- (0,1,0) node[anchor=south west] {$\CC \theta_1$};
                \draw[thick, ->] (0,0,0) -- (0,0,1) node[anchor=west] {$\CC \theta_2$};
            
                \fill[orange, opacity=0.3] (0,0,0) -- (1,0,0) -- (0,1,0) -- cycle; 
                \fill[orange, opacity=0.3] (0,0,0) -- (1,0,0) -- (0,0,1) -- cycle; 
                \fill[orange, opacity=0.3] (0,0,0) -- (0,1,0) -- (0,0,1) -- cycle; 
            
            \end{tikzpicture}
            \caption{The decorated fan of $X$}
        \end{subfigure}%
        \begin{subfigure}{.5\textwidth}
          \centering
          \begin{tikzpicture}[tdplot_main_coords, scale=2.5]
                \draw[->] (0,0,0) -- (1.5,0,0) node[anchor=north east] {$t_1$};
                \draw[->] (0,0,0) -- (0,1.5,0) node[anchor=north west] {$t_2$};
                \draw[->] (0,0,0) -- (0,0,1.5) node[anchor=south] {$t_3$};
            
                \draw[thick, ->] (0,0,0) -- (1,0,0) node[anchor=east] {$\CC(\theta_1-\theta_2)$};
                \draw[thick, ->] (0,0,0) -- (0,1,0) node[anchor=south west] {$\CC \theta_1$};
                \draw[thick, ->] (0,0,0) -- (0,0,1) node[anchor=west] {$\CC \theta_2$};
                \draw[thick, ->] (0,0,0) -- (0.8,0.8,0.8) node[anchor=west] {$0$};
            
                \fill[red, opacity=0.3] (0,0,0) -- (1,0,0) -- (0,1,0) -- (0.8,0.8,0.8) -- cycle; 
                \fill[blue, opacity=0.3] (0.8,0.8,0.8) -- (1,0,0) -- (0,0,1) -- cycle; 
                \fill[green, opacity=0.3] (0,0,0) -- (0,1,0) -- (0,0,1) -- cycle; 
            
            \end{tikzpicture}
          \caption{The decorated fan of $X'$}
          \label{}
        \end{subfigure}
        
        \caption{Decorated fans for a resolution of singularities $X' \to X$ as in Example \ref{ex:SingularButEvenSmooth}. The additional ray is generated by $(1,1,1)$ and decorated by the 0 subspace of $\t_\1$.}
        \label{fig:ResolutionOfSingularities}
    \end{figure}
\end{example}

\section{Orbit closures in the isomeric supergrassmannian}\label{QGrSection}

In this section we examine a natural collection of examples of toric supervarieties for quotients of the supertorus $T = Q(1)^n$, defined in Example \ref{Q(1)nExample}. We also illustrate a relationship between their decorated polytopes and some naturally-obtained subsets of $\t^*$.

\subsection{The isomeric supergrassmannian as a supervariety}

\begin{definition}
    An \textit{isomeric vector space} is a super vector space $V$ equipped with an odd involution $\Pi_V : V \to V$. Likewise, an \textit{isomeric subspace} of $V$ is a subspace which is preserved by $\Pi_V$.\footnote{We prefer the term \textit{isomeric}, as suggested in \cite{NSS}, to the more dated term \textit{queer}.}
\end{definition}

Write $Q(V,\Pi_V)$ for the supergroup of (even) automorphisms of $V$ that commute with $\Pi_V$, so $Q(V,\Pi_V) \cong Q(n)$ for $n = \dim V_\0 = \dim V_\1$.

\begin{definition}
    Let $(V,\Pi_V)$ be an isomeric vector space. The \textit{isomeric supergrassmannian} is the super moduli space $\QGr(r,V)$ of $(r|r)$-dimensional isomeric subspaces of $V$. Oftentimes we will assume $V = \CC^{n|n}$ with $\Pi_V : \CC^{n|0} \to \CC^{0|n}$ the ``identity" on $\CC^n$, in which case we will write $\QGr(r,n)$.
\end{definition}

For $A$ a superalgebra, the $A$-points of $\QGr(r,n)$ may be written as the space of full-rank matrices of the form
\begin{align*}
M = \begin{pNiceArray}{ccc|ccc}
    a_{11} & \cdots & a_{1r} & \alpha_{11} & \cdots & \alpha_{1r} \\
    \vdots & \ddots & \vdots & \vdots & \ddots & \vdots \\
    a_{n1} & \cdots & a_{nr} & \alpha_{n1} & \cdots & \alpha_{nr} \\
    \hline
    \alpha_{11} & \cdots & \alpha_{1r} & a_{11} & \cdots & a_{1r} \\
    \vdots & \ddots & \vdots & \vdots & \ddots & \vdots \\
    \alpha_{n1} & \cdots & \alpha_{nr} & a_{n1} & \cdots & a_{nr}
\end{pNiceArray},
\end{align*}
for $a_{ij} \in A_\0$ and $\alpha_{ij} \in A_\1$, modulo the right action of $Q(r)(A)$. As seen in \cite{Noja}, it is covered by $\binom{n}{r}$ affine charts (one for each size-$(n-r)$ subset $I$ of $\{1,...,n\}$) of the form
\begin{align*}
M' = \begin{pNiceArray}{c|c}
    a_I & \alpha_I \\
    1_{r \times r} & 0_{r \times r}  \\
    \hline
    \alpha_I & a_I \\
    0_{r \times r} & 1_{r \times r}
\end{pNiceArray},
\end{align*}
where $\begin{pNiceArray}{c|c} a_I & \alpha_I \end{pNiceArray}$ is the submatrix of $M$ consisting of the rows corresponding the indices in $I$, and where the rows of $M'$ have been rearranged accordingly for notational convenience. It is clear from this affine cover that $\QGr(r,n)$ is a smooth supervariety with underlying variety $\Gr(r,n)$.

\subsection{The isomeric supergrassmannian as a homogeneous space}
By way of the transitive action of $Q(n)$ on $\QGr(r,n)$, we see that $\QGr(r,n) \cong Q(n) / P$ for $P$ a maximal parabolic subgroup. The same holds if we replace $Q(n)$ by its subgroup $SQ(n)$ consisting of those linear transformations with $Q$-determinant 1, or even its quotient $PSQ(n)$ by scalars. It is from this construction that we view $\QGr(r,n)$ as the isomeric analog of the ordinary Grassmannian.

However, for the purposes of this paper, we are more interested in the action of the isomeric unitary supergroup $UQ(n)$ consisting of automorphisms which commute with $\Pi_V$ and preserve a super-Hermitian form $h$ in the sense of \cite{CCTV}, i.e.\ one that satisfies $h(u,v) = (-i)^{|u| \cdot |v|}\ol{h}(u,v)$ for a usual Hermitian form $\ol{h}$. From $h$ we obtain the notion of super-Hermitian adjoint
$$\begin{pmatrix}
    C & \Gamma \\ \Gamma & C
\end{pmatrix}^* = \begin{pmatrix}
    C^\dagger & -i \Gamma^\dagger \\ -i \Gamma^\dagger & C^\dagger
\end{pmatrix}$$
where $^\dagger$ represents the usual Hermitian adjoint with respect to $\ol{h}$.

We may therefore write the real Lie supergroup $UQ(n)$ as the one whose $A$-points (for $A$ a \textit{real} commutative superalgebra) are
$$UQ(n)(A) = \left\{ X  \in Q(n)(\CC \otimes_\RR A) \mid X^{-1} = X^* \right\}$$
and its Lie superalgebra
\begin{align*}
    \uq(n) &= \{X \in \q(n) \mid X + X^* = 0\} \\
    &= \left\{ \begin{pmatrix} C & \Gamma \\ \Gamma & C \end{pmatrix} \mid C = -C^\dagger, \Gamma = i\Gamma^\dagger \right\} \\
    &= \left\{ \begin{pmatrix} C & \Gamma \\ \Gamma & C \end{pmatrix} \mid C \in iH(n), \Gamma \in (1+i)H(n) \right\} \\
    &= iHQ(n)
\end{align*}
where $H(n)$ is the space of ordinary $n \times n$ Hermitian matrices, and
\begin{align*}
    HQ(n) &= \{X \in \q(n) \mid X^* = X\} \\
    &= \left\{ \begin{pmatrix} C & \Gamma \\ \Gamma & C \end{pmatrix} \mid C = C^\dagger, \Gamma = -i\Gamma^\dagger \right\} \\
    &= \left\{ \begin{pmatrix} C & \Gamma \\ \Gamma & C \end{pmatrix} \mid C \in H(n), \Gamma \in (1-i)H(n) \right\}
\end{align*}
is the space of isomeric super Hermitian endomorphisms. We may realize $\Pi HQ(n)$ (the parity shift of $HQ(n)$) as the dual space to $\uq(n)$ via the following odd pairing:
\begin{align*}
    \langle -,- \rangle : HQ(n) \times iHQ(n) &\to \RR \\
    \left\langle \begin{pmatrix} C & \Gamma \\ \Gamma & C \end{pmatrix}, \begin{pmatrix} C' & \Gamma' \\ \Gamma' & C' \end{pmatrix} \right\rangle &= \frac{1}{1+i} \tr(C\Gamma' + C'\Gamma)
\end{align*}

As in the ordinary situation, $UQ(n)$ acts transitively on $\QGr(r,n)$. The stabilizer of a fixed isomeric subspace is $UQ(r) \times UQ(n-r)$, so the following lemma holds:
\begin{lemma}\label{QGrUQnLemma}
    If $0 \leq r \leq n$, then $$\QGr(r,n) \cong UQ(n) / UQ(r) \times UQ(n-r).$$
\end{lemma}

\subsection{The isomeric supergrassmannian as an adjoint orbit}\label{sec:AdjointOrbitSuperpolytope}
Consider the adjoint action of $UQ(n)$ on $\uq(n) \cong iHQ(n)$:
\begin{align*}
    UQ(n)(A) \times \uq(n)(A) &\to \uq(n)(A) \\
    g \cdot X &\mapsto gXg^{-1}
\end{align*}
Here, $A$ is a real commutative superalgebra. Let us fix the closed point $$M_r = i\begin{pmatrix} \mathbbm{1}_r & 0 \\ 0 & \mathbbm{1}_r \end{pmatrix}$$ of $\uq(n)$ where $\mathbbm{1}_r$ is the diagonal matrix whose first $r$ diagonal entries are 1, with the rest 0. Since the stabilizer of $M_r$ is exactly $UQ(r) \times UQ(n-r)$, it follows by Lemma \ref{QGrUQnLemma} that the adjoint orbit $UQ(n) \cdot M_r$ can be equipped with a complex structure so that it is isomorphic to $\QGr(r,n)$. We will equivalently think of $\QGr(r,n)$ as an ``odd coadjoint orbit," i.e.\ the orbit of $-iM_r$ as a closed point in the parity shift $\Pi HQ(n)$.

We emphasize that there is no way to view $UQ(n) \cdot M_r$ as an even coadjoint orbit, since the adjoint and coadjoint representations of $UQ(n)$ are only isomorphic up to parity shift. In particular, the stabilizer of the closed point $$(1-i)\begin{pmatrix}
    0 & \mathbbm{1}_r \\
    \mathbbm{1}_r & 0
\end{pmatrix}$$ is not $UQ(r) \times UQ(n-r)$, but $U(r) \times UQ(n-r)$.

Denote by $(\t_\RR)_\0$ the Cartan subalgebra of $\uq_\0(n)$ consisting of diagonal matrices. Let $\t_\RR$ be its centralizer in $\uq(n)$, consisting of block matrices for which both $C$ and $\Gamma$ are diagonal. Observe that $\t_\RR$ is the Lie superalgebra of a compact real form of the supertorus $Q(1)^n \subset Q(n)$, and that $\t_\RR^*$ also consists of matrices for which both $C$ and $\Gamma$ are diagonal.

Consider the natural inclusion of supermanifolds $\QGr(r,n) \to \Pi \uq(n)^*$, which resembles an odd version of a momentum map. Composing it with the natural projection $\Pi \uq(n)^* \to \Pi \t_\RR^*$ yields a new map $\mu : \QGr(r,n) \to \Pi \t_\RR^*$. We now calculate the image of $\mu$ on $\RR[\xi]$-points so that we may understand the image of $\mu$ as a subset of the super vector space $\Pi \t_\RR^*$.

Write $C+\Gamma$ for $\left(\begin{smallmatrix}
    C & \Gamma \\ \Gamma & C
\end{smallmatrix}\right)$ in order to compress notation. Let $C+\Gamma \in UQ(n)(A)$, so $(C+\Gamma)^{-1} = C^\dagger - i\Gamma^\dagger$. Write $C = (c_{jk} + id_{jk})$ and $\Gamma = (\gamma_{jk} + i\delta_{jk})$. Then
\begin{align*}
    (C+\Gamma) \cdot (-iM_r) &= (C+\Gamma) (-iM_r) (C^\dagger - i\Gamma^\dagger) \\
    &= (c_{jk} + id_{jk} + \gamma_{jk} + i\delta_{jk})_{k \leq r} (c_{kj} - id_{kj} - \delta_{kj} - i\gamma_{kj}),
\end{align*}
a matrix whose $j$th diagonal entry is
\begin{align*}
    \sum_{k=1}^r \bigg( \big[c_{jk}^2 + d_{jk}^2 \big] + (1-i) \big[ -c_{jk}\delta_{jk} + d_{jk}\gamma_{jk} + c_{jk}\gamma_{jk} + d_{jk}\delta_{jk} \big] \bigg).
\end{align*}
Since $(C+\Gamma)(C^\dagger - i\Gamma^\dagger) = (C^\dagger - i\Gamma^\dagger)(C+\Gamma) = 1$, one may verify that
\begin{align*}
    \sum_{k=1}^n \bigg( \big[c_{jk}^2 + d_{jk}^2 \big] + (1-i) \big[ -c_{jk}\delta_{jk} + d_{jk}\gamma_{jk} + c_{jk}\gamma_{jk} + d_{jk}\delta_{jk} \big] \bigg) = 1
\end{align*}
for all $j$, and similarly
\begin{align*}
    \sum_{j=1}^n \bigg( \big[c_{jk}^2 + d_{jk}^2 \big] + (1-i) \big[ -c_{jk}\delta_{jk} + d_{jk}\gamma_{jk} + c_{jk}\gamma_{jk} + d_{jk}\delta_{jk} \big] \bigg) = 1
\end{align*}
for all $k$. Thus the sum of the diagonal entries of $(C+\Gamma) \cdot (-iM_r)$ is equal to $r$.

Write $d_1, ..., d_n$ for these diagonal entries. We find from the above equations that $\sum_{i=1}^n d_i = r$, and $d_i \in [0,1] + (1-i)\RR \xi$. Moreover, if the even part of such an entry is 0, then $c_{jk} = d_{jk} = 0$ for all $k=1, ..., r$, so the odd part is 0 as well. Likewise, if the even part is 1, then the odd part is the odd part of the trace, which is also 0. On the other hand, if the even part belongs to the open interval $(0,1)$, then the odd part can take on any value in $\RR \xi$.

In this sense, we may view the ``image" of $\mu : M \to \Pi\t_\RR^*$ as a decorated polytope $\P$ whose even part is the $(r,n)$-hypersimplex $\Delta(r,n)$. If $F$ is a face of $\Delta(r,n)$, let $I_F \subseteq \{1, ..., n\}$ be the collection of indices $i$ for which the $i$th coordinate of $F$ is 0 or 1. Then to the face $F$ we attach the subspace $W_F \subseteq \Pi (\t_\RR)_\0^* \cong \RR^n$ consisting of vectors whose entries sum to 0, and whose $i$th entry is 0 for any $i \in I$. Therefore
\begin{align*}
    \P &= \bigcup_{F \leq \Delta(r,n)} \Int F \times W_F \\
    &\cong \left\{(b_1, ..., b_n \mid \beta_1, ..., \beta_n) \in \RR^n \;\Bigg|\; (b_i, \beta_i) \in \{0,1\} \times \{0\} \cup (0,1) \times \RR \text{ and } \sum_{i=1}^n (b_i, \beta_i) = (r,0) \right\}
\end{align*}

\subsection{Supertorus actions}

By restricting the action of $Q(n)$ to the Cartan subgroup $Q(1)^n$ of diagonal matrices, we obtain a supertorus action on $\QGr(r,n)$. Let us consider its action (by matrix multiplication) on the $A$-points:

\begin{align*}
    \begin{pNiceArray}{ccc|ccc}
        b_1 &&& b_1 \beta_1 && \\
        & \ddots &&& \ddots & \\
        && b_n &&& b_n \beta_n \\
        \hline
        b_1 \beta_1 &&& b_1 && \\
        & \ddots &&& \ddots & \\
        && b_n \beta_n &&& b_n \\
    \end{pNiceArray} \cdot \begin{pNiceArray}{ccc|ccc}
    a_{11} & \cdots & a_{1r} & \alpha_{11} & \cdots & \alpha_{1r} \\
    \vdots & \ddots & \vdots & \vdots & \ddots & \vdots \\
    a_{n1} & \cdots & a_{nr} & \alpha_{n1} & \cdots & \alpha_{nr} \\
    \hline
    \alpha_{11} & \cdots & \alpha_{1r} & a_{11} & \cdots & a_{1r} \\
    \vdots & \ddots & \vdots & \vdots & \ddots & \vdots \\
    \alpha_{n1} & \cdots & \alpha_{nr} & a_{n1} & \cdots & a_{nr}
\end{pNiceArray},
\end{align*}
or equivalently
\begin{align*}
    \begin{pmatrix}
        b_1(1+\beta_1) &&\\
        & \ddots & \\
        && b_n(1+\beta_n)
    \end{pmatrix} \cdot \begin{pmatrix}
        a_{11}+\alpha_{11} & \cdots & a_{1r} + \alpha_{1r} \\
        \vdots & \ddots & \vdots \\
        a_{n1} + \alpha_{n1} & \cdots & a_{nr} + \alpha_{nr}
    \end{pmatrix}
\end{align*}
Although this action does not in general have an open orbit, we may still consider the closures of orbits of closed points as toric supervarieties for a quotient of $T = Q(1)^n$.

As before, we will write $\t = \q(1)^n = \CC\{x_1, ..., x_n\} \oplus \CC\{\theta_1, ..., \theta_n\}$ where $\frac{1}{2}[\theta_i, \theta_j] = \delta_{ij} x_i$. Abusing the language of Harish-Chandra pairs, we will often conflate a supertorus with its Lie superalgebra when writing stabilizer subgroups.

\begin{lemma}
    Let $y$ be a closed point of $\QGr(r,n)$.
    \begin{enumerate}
        \item $\Stab_{T_\0}(y)$ is a span of independent elements of the form $x_{i_1} + ... + x_{i_\ell}$, such that each $x_i$ occurs in exactly one such element.
        \item $\Stab_T(y)_\0 = \Stab_{T_\0}(y)$
        \item $\Stab_T(y)_\1 = \CC\{\theta_{i_1} + ... + \theta_{i_\ell} \mid x_{i_1} + ... + x_{i_\ell} \in \Stab_{T_\0}(y)\}$
    \end{enumerate}
\end{lemma}

\begin{proof}
    \begin{enumerate}
        \item Let $y$ be the closed point
        $$\begin{pmatrix}
            a_{11} & \cdots & a_{1r} \\
            \vdots & \ddots & \vdots \\
            a_{n1} & \cdots & a_{nr}
        \end{pmatrix} (\CC^\times)^r$$
        for $a_{ij} \in \kk$. Assume without loss of generality that $a_{ij} = \delta_{ij}$ for $i \leq r$, so that
        $$y = \begin{pmatrix}
            1 & & \\
            & \ddots & \\
            & & 1 \\
            a_{(r+1)1} & \cdots & a_{(r+1)r} \\
            \vdots & \ddots & \vdots \\
            a_{n1} & \cdots & a_{nr}
        \end{pmatrix} (\CC^\times)^r$$
        and hence $(t_1, ..., t_n) \in T_\0$ stabilizes $y$ if and only if $t_i = t_j$ whenever $a_{ij} \neq 0$.

        \item Immediate from the definitions.

        \item The same argument as in part (a) applies, using the compressed expression for the group action as matrix multiplication. \qedhere
    \end{enumerate}
\end{proof}

Write $\Pi_\t : \t \to \t$ for the odd involution $x_i \mapsto \theta_i$ induced by $\Pi_V$, so that $\Stab_T(y)$ is $\Pi_\t$-invariant. In particular, if $y$ belongs to the $T_\0$-orbit corresponding to $\sigma$, then $\Stab_T(y) = \t_{N_\sigma} \oplus \Pi \t_{N_\sigma}$.

\begin{proposition}
    Let $y$ be a closed point of $\QGr(r,n)$. The closure of the orbit $T \cdot y$ is a large-orbit toric supervariety corresponding to the decorated polytope $\P$ consisting of the usual even polytope $P$ of $\cl(T_\0 \cdot y)$ and odd subspaces $W_F = (\t_\1 / \Pi \t_{N_{\sigma_F}} )^*$, where $\sigma_F$ the normal cone to the face $F$ of $P$.
\end{proposition}

\begin{proof}
    Let $U$ be an open affine chart containing $y$; without loss of generality $$U = \left\{ \begin{pmatrix}
        1_{r \times r}  \\
        a_I
    \end{pmatrix} \right\}$$
    at the level of $\kk$-points. We have
    \begin{align*}
        &\begin{pmatrix}
            b_1(1+\beta_1) &&\\
            & \ddots & \\
            && b_n(1+\beta_n)
        \end{pmatrix} \cdot \begin{pmatrix}
            1 & & \\
            & \ddots & \\
            & & 1 \\
            a_{(r+1)1} & \cdots & a_{(r+1)r} \\
            \vdots & \ddots & \vdots \\
            a_{n1} & \cdots & a_{nr}
        \end{pmatrix} \\
        &= \begin{pmatrix}
            b_1(1+\beta_1) & & \\
            & \ddots & \\
            & & b_r(1+\beta_r) \\
            b_{r+1}(1+\beta_{r+1})a_{(r+1)1} & \cdots & b_{r+1}(1+\beta_{r+1}) a_{(r+1)r} \\
            \vdots & \ddots & \vdots \\
            b_n(1+\beta_n)a_{n1} & \cdots & b_n(1+\beta_n) a_{nr}
        \end{pmatrix} \\
        &= \begin{pmatrix}
            1 & & \\
            & \ddots & \\
            & & 1 \\
            b_1^{-1} b_{r+1} (1-\beta_1 + \beta_{r+1} + \beta_1\beta_{r+1}) a_{(r+1)1} & \cdots & b_r^{-1} b_{r+1} (1-\beta_r + \beta_{r+1} + \beta_r\beta_{r+1}) a_{(r+1)r} \\
            \vdots & \ddots & \vdots \\
            b_1^{-1} b_n (1-\beta_1 + \beta_n + \beta_1\beta_n) a_{n1} & \cdots & b_r^{-1} b_n (1-\beta_r + \beta_n + \beta_r\beta_n) a_{nr}
        \end{pmatrix}
    \end{align*}
    so that $\CC[U \cap \cl(T \cdot y)]$ is generated by $t_i^{-1} t_j (1+\xi_i \xi_j)$ and $t_i^{-1} t_j(-\xi_i+\xi_j)$ for $1 \leq i \leq r < j \leq n$. Hence $U \cap \cl(T \cdot y)$ is the large-orbit affine toric supervariety with $V_{\rho,0} = \Pi \t_{N_\rho}$ and $V_{\rho,1}=0$, so we are finished.
\end{proof}

Generically, up to complexification and composition with $\Pi_\t^* : \Pi \t^* \to \t^*$, this is the same decorated polytope obtained in section \ref{sec:AdjointOrbitSuperpolytope}! Other orbits correspond to decorated subpolytopes in a predictable manner. In effect, this means that the theory of torus strata developed in \cite{GGMS} generalizes easily to the situation of $Q(1)^n$ acting on $\QGr(r,n)$.

\appendix\section{Proof of Proposition/Definition \ref{propdef:HR1}}\label{section:ProofOfHR1}

In this section we prove equivalence of the following five conditions for a toric supervariety $X$ with decorated fan $(N,\Sigma,\t,\h,\{V_{\rho,i}\})$.
\begin{enumerate}
    \item For every $\rho \in \Sigma(1)$, there is a parameterization as in (\ref{eqn:LeftDerivations}) of the derivations by which $\t$ acts on $\CC[T]$ such that
    $$A_\rho = \CC[t_1, t_2^{\pm 1}, ..., t_p^{\pm 1}, \xi_1, ..., \xi_r, t_1^{\ell_{r+1}} \xi_{r+1}, ..., t_1^{\ell_q} \xi_q]$$
    for some positive integers $\ell_i$.
    
    \item For every $\rho \in \Sigma(1)$, $A_\rho$ is the subalgebra of $\CC[T]$ generated by $L_{\CC[T]^W}(n+\rho^\perp)$ for all $n \geq 0$ and all codimension-at-most-1 subspaces $W \subseteq V_{\rho,0}$ containing $V_{\rho,n}$.

    \item For every $\rho \in \Sigma(1)$ and proper subspace $W \subset V_{\rho,0}$ containing $\h$, the quotient $T/W$ of $T/H$ embeds into a toric supervariety with decorations $(N,\rho,\t,W,\{V_{\rho,i}+W\})$.

    \item For every point $x \in |X|$ of codimension at most 1, the maximal ideal $\m_x \subseteq \O_{X,x}$ is generated by a regular sequence $r_1, ..., r_d$ such that whenever $\theta \in (\Stab_\t x)_\1$ and $\ol{\theta \cdot r_i} = 0$, we have $\theta \cdot r_i \in (r_i)$.

    \item For every point $x \in |X|$ of codimension at most 1 and for every $\theta \in (\Stab_\t x)_\1$, $DS_\theta (\O_{X,x}/((\O_{X,x}^{\t_\0})_\1))$ is FR.
\end{enumerate}

First, we need the following lemma. In essence, it says that locally in codimension 1, the badness of a non-$(HR_1)$ toric supervariety is contained in the badness of the case of an abelian supertorus. That is, each $\Spec A_\rho$ contains a closed subvariety which is a toric supervariety for the abelian supertorus $(T_\0, \t_\0 \oplus V_{\rho,0})$. This closed toric subvariety is the subject of condition (e) of Proposition/Definition \ref{propdef:HR1}.

\begin{lemma}\label{lemma:RayStructureAppendix}
Let $\Spec A_\rho$ be a toric supervariety with $[V_{\rho,0}, V_{\rho,0}]=0$, and parameterize the derivations so that $\theta_j^\r = \d{\xi_j}$ for $\theta_j \in V_{\rho,0}$. Then we may write
$$A_\rho = \CC[t_1, t_2^{\pm 1}, ..., t_p^{\pm 1}, \xi_1, ..., \xi_r, t_1^{\ell_{r+1}} (\xi_{r+1}+...), ..., t_1^{\ell_{q}} (\xi_{q}+...)]$$
such that each $(\xi_j+...)$ for $j = r+1, ..., q$ lies in $\bigwedge\{\xi_{r+1}, ..., \xi_q\}$.
\end{lemma}

\begin{proof}
By default, we know
$$A_\rho = \CC[t_1(1+...), t_2^{\pm 1}, ..., t_p^{\pm 1}, \xi_1, ..., \xi_r, t_1^{\ell_{r+1}} (\xi_{r+1}+...), ..., t_1^{\ell_{q}} (\xi_{q}+...)]$$
without any assumptions on the ellipses. Our assumption on $\theta_j^\r$ is equivalent to specifying $x_{ij}=0$ for $j=r+1, ..., q$ in the language of \ref{eqn:LeftDerivations}. It follows that for $j=1, ..., q$, we have
\begin{align*}
    \theta_j \cdot t^m(...) \in \CC t^m(1+...) \bigwedge\{\xi_1, ..., \xi_r\} - \d{\xi_j} t^m(...),
\end{align*}
of which the first summand is contained in $A_\rho$. Hence the second summand is too, so in particular $A_\rho$ is closed under the action by $\d{\xi_j}$. 

By repeated application of these derivations, we conclude that $t_1 \in A_\rho$. Likewise, if some $t_1^{\ell_j}(\xi_j+...)$ contains a summand with a $\xi_i$ term for $i=1, ..., r$, then we may subtract $\xi_i \d{\xi_i} t_1^{\ell}(\xi_j+...)$ from it so that the result does not.
\end{proof}

We are now ready to prove Proposition/Definition \ref{propdef:HR1}.

\begin{proof}
\textbf{Case 1, (a)$\iff$(b)$\iff$(c):} First suppose $[V_{\rho,0}, V_{\rho,0}] \neq 0$, so $\codim(\h, V_{\rho,0}) = 1$ and $V_{\rho,1} = \h$. We claim that conditions (a)-(e) always hold in this case. First note that we may pick a basis $\theta_1, ..., \theta_s$ of $\t_\1$ with dual basis $\xi_1, ..., \xi_s$ so that $\h = \CC\{\theta_{q+1}, ...,\theta_s\}$ and $V_{\rho,0} = \CC\{\theta_q, ..., \theta_s\}$. We may then parameterize $\CC[T]$ so that
\begin{align*}
    \theta_q^\r &= \sum_{i=q}^s a_i \xi_i t_1 \d{t_1} + \d{\xi_q} & \\
    \theta_j^\r &= \d{\xi_j} & j=q+1, ..., s
\end{align*}
for some $a_i \in \CC$. Therefore $L_{A_\rho}(\rho^\perp) = \CC[T/H_{\rho,0}] = \CC[t_2^{\pm 1}, ..., t_p^{\pm 1}, \xi_1, ..., \xi_{q-1}]$, and there is a new odd generator in $L(1)$. Since $L(1)$ must be $2^q$-dimensional, it follows that $L_{A_\rho}(1+\rho^\perp) = L_{\CC[T/H]}(1+\rho^\perp)$ and hence both (a) and (b) follow. In particular, we have
$$A_\rho = \CC[t_1, t_2^{\pm 1}, ..., t_p^{\pm 1}, \xi_1, ..., \xi_{q-1}, t_1 \xi_q].$$ Note that (c) is trivial since the only possible $W$ is $\h$ itself.

\textbf{A brief aside:} We show that for the purposes of checking parts (d) and (e), even in the case of $[V_{\rho,0}, V_{\rho,0}]=0$, it suffices to consider $T_\0$-invariant points $x \in |X|$ of codimension 1. Such a point which is not $T_\0$-invariant can be viewed as a prime ideal (of some $\ol{A_\sigma}$) generated by a single (regular) element which is not a monomial. As a result, the local ring $\O_{X,x}$ will contain all $\xi_1, ..., \xi_q$, so $(\Stab_\t x)_\1 = \h$. This conclusion also holds for the codimension-0 point.

Since $\theta_j^\r = \d{\xi_j}$ kills $\O_{X,x}$ for $\theta_j \in \h$, we note that $\theta_j$ simply takes each $t^m$ to some multiple of it by a linear combination of $\xi_1, ..., \xi_q$. Hence (d) will always hold for this choice of $x$. Meanwhile, (e) is trivial since $\O_{X,x}^{\t_\0}$ contains all $\xi_1, ..., \xi_q$.

On the other hand, a $T_\0$-invariant point of codimension 1 can be viewed as the ideal generated by $t_1$ and the odd elements in some $A_\rho$. Since localization at this prime will not affect future calculations, we will often suppress it from the notation and simply calculate with $A_\rho$.

\textbf{Case 1, (a)$\iff$(d)$\iff$(e):}
Let us now prove (d) in this case. We have $(\Stab_\t x)_\1 = \CC \theta_q + \h$. By the same argument as in the aside, any contribution of $\h$ here is insignificant. Hence we check that the condition holds for $\theta=\theta_q$ and $r_1, ..., r_d = t_1, \xi_1, ..., \xi_{q-1}, t_1 \xi_q$. We have $\theta_q \cdot t_1 = \sum_{i=1}^{q} b_i \xi_i t_1$ for some $b_i \in \CC$ and $\theta_q \cdot \xi_i = 0$ for all $i=1, ..., q-1$, as well as $\theta_q \cdot t_1 \xi_q = t_1+...$, so (d) holds. For (e), we note that the ideal $((\O_{X,x}^{\t_\0})_\1)$ is generated by $\xi_1, ..., \xi_{q-1}$, and $A_\rho^{\theta_q^2} = \CC[t_2^{\pm 1}, ..., t_p^{\pm 1}, \xi_1, ..., \xi_{q-1}]$. Thus the quotient is purely even, so we are finished.

\textbf{Case 2, (a)$\iff$(b)$\iff$(c):} We now assume $[V_{\rho,0}, V_{\rho,0}]=0$. Suppose (a) holds, so we may write $\theta_j^\r = \d{\xi_j}$ for $j=r+1, ..., s$, where $\h = \CC\{\theta_{q+1}, ..., \theta_s\}$. By default, we have $L_{A_\rho}(\rho^\perp) = \CC[T/H_{\rho,0}] = \CC[t_2^{\pm 1}, ..., t_p^{\pm 1}, \xi_1, ..., \xi_r]$. Now for $W = \CC\{\theta_{r+1}, ..., \hat{\theta_j}, ..., \theta_q\}$, we obtain $L_{\CC[T]^W} = \CC[t_1^{\pm 1}, ..., t_p^{\pm 1}, \xi_1, ..., \xi_r, \xi_j]$ and so (b) follows.

Conversely, if (b) holds, choose a basis $\xi_1, ..., \xi_s$ of $\t^*$ compatible with the chain $V_{\rho,0} \supseteq V_{\rho,1} \supseteq ...$ in the sense that each $V_{\rho,i}$ is spanned by a subset of the dual basis. In particular, assume $\xi_1, ..., \xi_r$ vanish on $V_{\rho,0}$, and $\xi_1, ..., \xi_q$ vanish on $\h$. Likewise, choose even coordinates $t_1, t_2^{\pm 1}, ..., t_q^{\pm 1} \in \CC[T]^{V_{\rho,0}}$, so that $\theta_i^{\r} = \d{\xi_i}$ for $i=r+1, ..., q$ in terms of these elements. If $W \supseteq V_{\rho,n}$ is a codimension-1 subspace of $V_{\rho,0}$, then $\CC[T]^W = \CC[s_1^{\pm 1}, ..., s_p^{\pm 1}, \xi_1, ..., \xi_r, \xi_k]$ for $\xi_k$ a functional vanishing on $W$ but not $V_{\rho,0}$. It follows that $L_{\CC[T]^W}(n+\rho^\perp)$ introduces the generator $s^n \xi_k$, so (a) is satisfied.

The equivalence of (b) and (c) is immediate.

\textbf{Case 2, (a)$\iff$(d)$\iff$(e):}
First assume (a), so
$$A_\rho = \CC[t_1, t_2^{\pm 1}, ..., t_p^{\pm 1}, \xi_1, ..., \xi_r, t_1^{\ell_{r+1}} \xi_{r+1}, ..., t_1^{\ell_q} \xi_q]$$
and
\begin{align*}
    \theta_j^\r &= \d{\xi_j} \\
    \theta_j \cdot t^m &= \left( \sum_{i=1}^r \langle m, [\theta_i,\theta_j] \rangle \xi_i \right) t^m \\
    \theta_j \cdot \xi_i &= -\delta_{ij}
\end{align*}
for $j=r+1, ..., s$. Now consider the regular sequence $t_1, \xi_1, ..., \xi_r, t_1^{\ell_{r+1}} \xi_{r+1}, ..., t_1^{\ell_q} \xi_q$. We have
\begin{align*}
    \theta_j \cdot t_1 &\in (t_1) \\
    \theta_j \cdot \xi_i &= 0 &i=1, ..., r \\
    \theta_j \cdot t_1^{\ell_{i}} \xi_i &\in (t_1^{\ell_i} \xi_i) &j \neq i >r \\
    \ol{\theta_j \cdot t_1^{\ell_j} \xi_j} &\neq 0
\end{align*}
for $j>r$, so indeed (d) holds. Likewise, we have $((\O_{X,x}^{\t_\0})_\1) = (\xi_1, ..., \xi_r)$, so checking (e) reduces to calculating the cohomology of $\d{\xi_j}$ on $\CC[t_1, t_2^{\pm 1}, ..., t_p^{\pm 1}, t_1^{\ell_{r+1}} \xi_{r+1}, ..., t_1^{\ell_q} \xi_q]$. This is exactly $\CC[t_1, t_2^{\pm 1}, ..., t_p^{\pm 1}, t_1^{\ell_{r+1}} \xi_{r+1}, ..., \widehat{t_1^{\ell_{j}} \xi_{j}},..., t_1^{\ell_q} \xi_q]/(t_1^{\ell_j})$, which is certainly FR.

On the other hand, suppose (a) is not true, and write
$$A_\rho = \CC[t_1, t_2^{\pm 1}, ..., t_p^{\pm 1}, \xi_1, ..., \xi_r, t_1^{\ell_{r+1}} (\xi_{r+1}+...), ..., t_1^{\ell_{q}} (\xi_{q}+...)]$$
as in Lemma \ref{lemma:RayStructureAppendix}. Assume without loss of generality that modulo $((A_\rho)_\1^2)$, the regular sequence $r_1, ..., r_d$ is a permutation of the generators of $A_\rho$ above. Now choose $\ell_k$ minimal such that the corresponding generator $t_1^{\ell_k}(\xi_k+...)$ contains a summand of the form $\xi_J$ for $\sum_{j \in J} \ell_j > \ell_k$. Then for $j \in J \subseteq \{r+1, ..., q\} \backslash \{k\}$, we have
\begin{align*}
    \theta_j \cdot t_1^{\ell_k}(\xi_k+...) &= \left( \sum_{i=1}^r \langle \ell_k, [\theta_i, \theta_j] \rangle \xi_i \right)  t_1^{\ell_k}(\xi_k+...) - t_1^{\ell_k} \d{\xi_j}(...)
\end{align*}
wherein the first summand belongs to $(t_1^{\ell_j}(\xi_j+...))$ but the second does not. Hence (d) fails in this situation.

Finally, to see that (e) does not hold, we may reduce to the case in which $r=0$ and $\theta_j = \d{\xi_j}$. As above, let $\ell_k, J$, and $j$ be as above. Then the image of $\d{\xi_j}$ contains $t_1^{\ell_k}(\xi_{J \backslash j} + ...)$, meaning there is a relation in $DS_\theta \O_{X,x}$ among prior generators $t_1^{\ell_i}\xi_i$ for $i \in J \backslash \{j\}$. Hence the result is not FR, so we are finished.
\end{proof}

\bibliographystyle{hsiam}
\bibliography{refs}

\end{document}